\documentclass[12pt]{article}
\usepackage{amssymb}
\usepackage{amsmath}
\usepackage{array}

\usepackage{minitoc}

\begin{document}
\input{amssym.def}

\newsymbol \circledarrowleft 1309

\newcommand{\ha}{\mathfrak h}
\newcommand{\g}{\mathfrak g}
\newcommand{\ta}{\mathfrak t}
\newcommand{\s}{\mathfrak s}
\newcommand{\ctext}[1]{\makebox(0,0){#1}}
\setlength{\unitlength}{0.1mm}

\newcommand{\wt}{\widetilde}
\newcommand{\Lr}{\Longrightarrow}
\newcommand{\Aut}{\mbox{{\rm Aut}$\,$}}
\newcommand{\ul}{\underline}
\newcommand{\ol}{\overline}
\newcommand{\lr}{\longrightarrow}
\newcommand{\bc}{{\Bbb C}}
\newcommand{\bp}{{\Bbb P}}
\newcommand{\cf}{{\cal F}}

\newcommand{\ce}{{\cal E}}
\newcommand{\co}{{\cal O}}
\newcommand{\cg}{{\cal G}}
\newcommand{\hra}{\hookrightarrow}

\newtheorem{guess}{Theorem}[section]
\newcommand{\bth}{\begin{guess}$\!\!\!${\bf .}~~\rm}
\newcommand{\eeth}{\end{guess}}

\newtheorem{propo}[guess]{Proposition}
\newcommand{\bprop}{\begin{propo}$\!\!\!${\bf .}~~\rm}
\newcommand{\eprop}{\end{propo}}

\newtheorem{rema}[guess]{\it Remark}
\newcommand{\brem}{\begin{rema}$\!\!\!${\it .}~~\rm}
\newcommand{\erem}{\end{rema}}
\newtheorem{coro}[guess]{Corollary}
\newcommand{\bcor}{\begin{coro}$\!\!\!${\bf .}~~\rm}
\newcommand{\ecor}{\end{coro}}
\newtheorem{lema}[guess]{Lemma}
\newcommand{\blem}{\begin{lema}$\!\!\!${\bf .}~~\rm}
\newcommand{\elem}{\end{lema}}
\newtheorem{exam}{Example}%[section]
\newcommand{\beg}{\begin{exam}$\!\!\!${\bf .}~~\rm}
\newcommand{\eeg}{\end{exam}}
\newcommand{\er}{\hfill {\Large $\bullet$}\linebreak}
\newtheorem{defe}[guess]{Definition}
\newcommand{\bdefe}{\begin{defe}$\!\!\!${\bf .}~~ \rm}
\newcommand{\edefe}{\end{defe}}
\newcommand{\spec}{{\rm Spec}\,}

\title{Semistable Principal Bundles-II\\ (in positive
  characteristics)} \author{V.Balaji\footnote {The research of the
    first author was partially supported by the DST project no
    DST/MS/I-73/97}~ and A.J.Parameswaran} \date{} \maketitle

\vspace{2mm}
\noindent

\section{Introduction}
Let $H$ be a semisimple algebraic group and let $X$ be a smooth
projective curve defined over an algebraically closed field $k$.

One of the important problems in the theory of principal $H$-bundles
on $X$ is the construction of the moduli spaces of semistable
$H$-bundles when the characteristic of $k$ is positive. Over fields of
characteristic $0$ this work was done by A.Ramanathan (cf.\cite{r1}).
For principal $GL(n)$-bundles this is classical, over fields of any
characteristic (cf.\cite{ses}).

The purpose of this paper is to prove the existence and the
projectivity of the moduli spaces of semistable principal $H$-bundles
on $X$ for fields $k$ of characteristic $p > 0$ with precise bounds on
the prime $p$, the restrictions being imposed by the representation
theory of $H$.

It might seem, by the general method of reduction modulo $p$, that the
existence of the moduli space in char.0 implies its existence for
large primes. To the best of our knowledge this is not the case. (cf
Remark \ref{larhep}). The representation theoretic considerations
involving heights are essential to the proving of the existence of the
moduli.

The broad strategy of this paper is along the same lines as in the
precursor to this paper (\cite{bals}) where a different approach for
the construction and projectivity of these moduli spaces (in
characteristic zero) was given. However, its implementation involves
several new inputs.  The key input for the {\it existence} of the
moduli comes from the paper of Ilangovan-Mehta-Parameswaran
(\cite{imp}) which establishes in positive characteristics the links
between the semistability of principal bundles and the concept of a
{\it low height representation}. In proving the {\it projectivity} of
the moduli space, the key ideas come from a natural interplay of the
recent results of Serre on the representation theory in positive
characteristics (\cite{mour}, \cite{tens}), and ideas inspired by the
papers of Ramanan-Ramanathan and Rousseau (\cite{rr}, \cite{rou}).
The principal difficulty is to replace {\it the tensor product theorem
  of semistable bundles} and {\it unitary representations of
  fundamental groups} which are central to the characteristic 0
theory. The notions of height and saturated groups provide just the
right replacements.

Let $H$ be a semisimple algebraic group (as coming by reduction from a
Chevalley group scheme defined over $\bf Z$), and fix a faithful
representation $H \hookrightarrow G = SL(n)$ arising as reduction
modulo $p$ of a representation defined over $\bf Z$. Let us denote by
$ht_H(G)$ the height of $G$ as an $H$-representation.  (cf.
Definition \ref{height}). We say a representation $H \hookrightarrow
G$ is of {\it low height} if $char(k) = p > ht_H(G)$. Then we
have the following:

{\bf Theorem \ref{moduli}} Let $H \hookrightarrow G$ be a faithful low
height representation. Then there exists a coarse moduli scheme
$M_X(H)$, for semistable principal $H$-bundles. Further, the moduli
space $M_X(H)$ is quasi-projective and the canonical morphism
$\ol{\mu}: M_X(H) \lr M_X(G)$ is {\it affine}.

The proof of the projectivity of the moduli spaces requires more
refined prime bounds. Towards this we introduce a new index which we
term the {\it separable index} associated to a $G$-module $W$ (cf.
Definition \ref{sepind}).  We denote this by $\psi_G(W)$ and we say a
$G$-module $W$ is of {\it low separable index} if $char(k) = p >
\psi_G(W)$. We fix throughout, a finite dimensional $G$-module $W$
such that the subgroup $H$ is realised as the isotropy of a closed
orbit, hence giving rise to a closed embedding $G/H \subset W$. We
term these $G$-modules for convenience as {\em affine $(G,H)$-modules}
(cf. Def \ref{chevmod}). Let $A$ be a complete discrete valuation ring
and let $K$ be its quotient field and $k$ its residue field. Then we
have the following theorem:

{\bf Theorem \ref{ssred}} ({\em Semistable reduction}) Let $W$ be a
finite dimensional {\em affine $(G,H)$-module} associated to $H$ and $G$
and let $p > \psi_G(W)$.  Let $H_K$ denote the group scheme $H \times
\spec K$, and $P_K$ be a semistable $H_K$-bundle on $X_K$.  Then there
exists a finite extension $L/K$, with $B$ as the integral closure of
$A$ in $L$ such that the bundle $P_K$, after base change to $\spec B$,
extends to a semistable $H_B$-bundle $P_B$ on $X_B$.

This in particular implies that the moduli spaces $M_X(H)$ are
projective over fields $k$ with $char(k) = p > \psi_G(W)$. Together
with Theorem \ref{moduli} we can conclude that the canonical morphism
$\ol{\mu}: M_X(H) \lr M_X(G)$ is {\it finite}. As a corollary we also
obtain the irreducibility of the moduli spaces when $H$ is semisimple
and simply connected.(Cor \ref{irred}). A large part of this paper is
devoted to proving Theorem \ref{ssred}.

The crucial difference between the present approach and the classical
proof of Langton for the properness of the moduli space of semistable
vector bundles can be briefly described as follows. Langton first
extends the family of semistable vector bundles (or equivalently
principal $GL(n)$-bundles) to a $GL(n)$-bundle in the limit although
non-semistable. In other words, the structure group of the limiting
bundle remains $GL(n)$. Then by a sequence of {\it Hecke modifications}
the semistable limit is attained without changing the isomorphism
class of the bundle over the generic fibre.

Instead, we extend the family of semistable $H_K$-bundles to an
$H_A'$-bundle with the limiting bundle remaining semistable, but the
structure group scheme $H_A'$, is non-reductive in the limit. In other
words one loses the reductivity of the structure group scheme. Then,
by using Bruhat-Tits theory (cf \S10), we relate the group scheme
$H_A'$ to the reductive group scheme $H_A$ without changing the
isomorphism class of the bundle over the generic fibre as well as the
semistability of the limiting bundle.

We note that the boundedness of semistable principal bundles over
curves in positive characteristics is proved in the preprint
(\cite{hn}).\footnote{The problem of the construction of the moduli is
being considered independently by V.B.Mehta and S.Subramaniam.}

Throughout the paper, we make an effort to specify carefully the
bounds on the characteristic of $k$ that are forced on us. We believe
that our methods can probably be stretched to include more primes and
we indicate at every stage the possible difficulties. The
representation theoretic indices that we have developed here may
possibly be of independent interest.

Before we proceed to describe the contents we pause to remark
that there is some overlap between the present paper and \cite{bals}.

The layout of the paper is as follows. In \S3 we recall low height
representations and some results from \cite{mour} which we need in later
sections. Here we also define the basic functors for semistable principal
$H$- and $G$-bundles and we prove a technical lemma involving the choice of
a ``base point on the curve''which, in some sense gives the motivation for
the rest of the work. In this paper we work with more than one base point
so as to achieve better height bounds.

In \S4 we give a simple construction of the moduli space of $H$-bundles
under the right characteristic bounds. The idea of the proof comes from
\cite{bals} and the ingredients involving heights from \S3.

The rest of the paper is devoted to proving the semistable reduction
theorem. In \S5 some new representation theoretic indices are
introduced and these give the bounds that we need to impose on the
characteristic $p$ in what follows. Here the main point is to give a
criterion for the {\it strong separability} of a linear action of a
reductive group. In sections \S6 and \S7 we construct and study the
flat closure of $H_K'$ in $G_A$ and realise it as {\it isotropy group
  schemes} along the lines of the classical theorem on
semi-invariants (cf. \cite{bore}).

In \S8 we prove the key lemmas on the relationship between polystable
bundles and semistable sections inspired mainly by the papers of
Ramanan-Ramanathan (\cite{rr}) and Rousseau (\cite{rou}). More precisely,
we obtain a notion paralleling that of {\it monodromy} subgroup of a
polystable $G$-bundle which is realised as a {\it saturated subgroup} of
$G$. This enables us to prove a {\it local constancy} for polystable
bundles in char.$p$. In \S9 we prove that the family of bundles extends to
a semistable bundle with structure group as a non-reductive group scheme
$H_A'$ with generic fibre $H_K$. In \S10 using Bruhat-Tits theory we relate
the non-reductive group scheme $H_A'$ with the reductive group scheme
$H_A$. In \S11 we complete the proof of the semistable reduction theorem.

\footnotesize{\it Acknowledgments}: We would like to thank the many
people with whom we have had discussions during the course of this
work: S.Kannan, M.S.Narasimhan, M.V.Nori, Gopal Prasad,
M.S.Raghunathan, C.S.Rajan, S.Ramanan, S.Ilangovan, S.Subramaniam, and
V.Uma.  We want to especially thank V.B.Mehta and C.S.Seshadri for
their generous help in the paper, from its inception to its
conclusion. Finally we wish to thank the referees for their numerous
suggestions and comments which has led to a considerable improvement
in the exposition.

The first author wishes to thank the School of Mathematics T.I.F.R,
Mumbai and the second author the Chennai Mathemaical Institute and
Institute of Mathematical Sciences, Chennai where much of this work
was carried out. We also wish to thank CAAG for its annual meets where
we got together on this project.
%\pagebreak
\footnotesize 
%\tableofcontents
\linebreak

{\large {\bf Contents}}
\begin{enumerate}
\item Introduction
\item Notations and Conventions
\item Low height representations 
\item Construction of moduli
\item Separable index and slice theorem
\item Towards the flat closure
\item Affine embedding of $G_A/H_A'$
\item Semistable bundles, semistable sections and saturated groups
\item Extension to the flat closure
\item Potential good reduction
\item Semistable reduction theorem
\end{enumerate}

\normalsize
\section{Notations and Conventions}
Throughout this paper, unless otherwise stated, we have the following
notations and assumptions:

{\renewcommand{\labelenumi}{{\rm (\roman{enumi})}}
\begin{enumerate}
\item We work over an algebraically closed field $k$ of 
characteristic $p > 0$.

\item $H$ is a {\it semisimple} algebraic group, and $G$, unless
otherwise stated will always stand for the special linear group
$SL(n)$.  Their representations are finite dimensional and rational.

\item $A$ is a discrete valuation ring (which could be assumed to be
complete) with residue field $k$, and quotient field $K$.

\item We recall that $\pi:E \lr X$ is a principal bundle with structure
group $H$, or a principal $H$-bundle for short if $H$ acts on $E$ on the
right and $\pi$ is $H$-invariant and isotrivial, i.e, locally trivial in
the \'etale topology.

\item Let $E$ be a principal $G$-bundle on $X \times T$ where $T$ is
  $\spec A$. Let $x \in X$ be a closed point which we fix throughout
  and we shall denote by $E_{x,A}$ or $E_{x,T}$ (resp $E_{x,K}$) the
  restriction of $E$ to the subscheme $x \times~\spec A$ or $x \times
  T$ (resp $x \times~\spec K$). Similarly, $l \in T$ will denote the
  closed point of $T$ and the restriction of $E$ to $X \times l$ will
  be denoted by $E_l$.

\item We shall denote $T-l$ by $T^*$ throughout this paper. 

\item In the case where the structure group is $GL(n)$, when we speak of a
principal $GL(n)$-bundle we identify it often with the associated vector
bundle (and can therefore talk of the degree of the principal
$GL(n)$-bundle).

\item We denote by $E_K$ (resp $E_A$) the principal bundle $E$ on $X
\times \spec K$ (resp $X \times \spec A$) when viewed as a principal
$H_K$-bundle (resp $H_A$-bundle). Here $H_K$ and $G_K$ (resp $H_A$ and
$G_A$) are the product group scheme $H \times \spec K$ and $G \times
\spec K$ (resp $H \times \spec A$ and  $G \times \spec A$).

\item If $H_A$ is an $A$-group scheme, then by $H_A(A)$ (resp
$H_K(K)$) we mean its $A$ (resp $K$)-valued points. When $H_A = H
\times \spec A$, then we simply write $H(A)$ for its $A$-valued
points. We denote the closed fibre of the group scheme by $H_k$.

\item Let $Y$ be any $G$-scheme and let $E$ be a $G$-principal
bundle. For example $Y$ could be a $G$-module. Then we denote by
$E(Y)$ the associated bundle with fibre type $Y$ which is the
following object: $E(Y)$ = $(E \times Y)/G$ for the twisted action of
$G$ on $E \times Y$ given by $g.(e,y)~=~(e.g,g^{-1}.y)$.

\item If we have a group scheme $H_A$ (resp $H_K$) over
$\spec A$ (resp $\spec K$) an $H_A$-module $Y_A$ and a
principal $H_A$-bundle $E_A$, then we shall denote the
associated bundle with fibre type $Y_A$ by $E_A(Y_A)$.

\item By a family of $H$ bundles on $X$ parametrised by $T$ we mean a
principal $H$-bundle on $X \times T$, which we also denote by
$\{E_t\}_{t \in T}$.
\end{enumerate}

\section{Low height representations and some consequences}

Let $k$ be an algebraically closed field of characteristic $p>0$. Let
$H$ be a connected reductive algebraic group over $k$.  Let $T$ be a
maximal torus of $H$, $X(T):=Hom(T,{\bf G}_m)$ be the character group of $T$
and $Y(T):=Hom({\bf G}_m,T)$ be the 1-parameter subgroups of $T$. Let
$R\subset X(T)$ be the root system of $H$ with respect to $T$. Let
${\cal W}$ be the Weyl group of the root system $R$.  Let $(~,~)$ denote
the ${\cal W}$-invariant inner product on $X(T)\otimes {\bf R}$. For 
$\alpha\in R$, the corresponding co-root $\alpha^{\vee}$ is
$2\alpha/(\alpha,\alpha)$.  Let $R^{\vee}\subset{X(T)\otimes {\bf
R}}$ be the set of all co-roots. Let $B\subset H$ be a Borel subgroup
containing $T$. This choice defines a base $\Delta^+$ of $R$ called
the {\it simple roots}. Let $\Delta^-=-\Delta^+$. A root in $R$ is
said to be {\em positive} if it is a non-negative linear combination
of simple roots. We take the roots of $B$ to be positive by
convention. Let $\Delta^{\vee}\subset R^{\vee}$ be the basis for 
the corresponding dual root system. Then we can define the Bruhat
ordering on ${\cal W}$. The longest element with respect to this
ordering of ${\cal W}$ is denoted by $w_0$.  A reductive group is
classified by these {\it root-data}, namely the character group,
1-parameter subgroups, the root system, co-roots and the ${\cal
W}$-invariant pairing.

Let $V$ be a $H$-module, i.e., $V$ is a $k$-vector space together with a
linear representation of $H$ in Aut $(V)$. Then $V$ can be written as
direct sum of eigenspaces for $T$. On each eigenspace $T$ acts by a
character. These are called the {\em weights} of the representation. A
weight $\lambda$ is called {\em dominant} if
$(\lambda~,~\alpha_i^{\vee})\geq 0$ for all simple roots $\alpha_i\in
\Delta^+$. A weight $\lambda$ is said to be $\geq$ another weight $\mu$ if
the difference $\lambda-\mu$ is a non-negative integral linear combination
of simple roots, where the difference is taken with respect to the natural
abelian group structure of $X(T)$. The {\em fundamental weights} $\omega_i$
are uniquely defined by the criterion
$(\omega_i~,~\alpha_j^{\vee})=\delta_{ij}$. The element $\rho$ of
$X(T)\otimes {\bf R}$ is defined to be half the sum of positive roots.  It
can also seen to be equal to the sum of fundamental weights.  The {\em
height} (cf. [H], Section 10.1) of a root is defined to be the sum of the
coefficients in the expression $\alpha~=~\Sigma k_i\alpha_i$. We extend
this notion of height linearly to the weight space and denote this function
by $ht(~)$. Note that $ht$ is defined for all weights but need not be an
integer even for dominant weights. We extend this notion of height to
representations as follows:
\linebreak
\bdefe\label{height} 
\begin{enumerate}
\item Given a linear representation $V$ of $H$, we define the
{\bf height} of the representation $ht_H(V)$ (also denoted by $ht(V)$ if
$H$ is understood in the given context) to be the maximum of
$2ht(\lambda)$, where $\lambda$ runs over dominant weights
occurring in $V$.  

\item A linear representation $V$ of $H$ is said to be a {\bf low
    height} representation if $ht_H(V)<p$, and a weight $\lambda$ is of
  low height if $2ht(\lambda)<p$.  
\end{enumerate}
\edefe

Then we have the following theorem (cf. \cite{imp}, \cite{mour})

\bth
Let $V$ be a linear representation of $H$ of low height. Then
$V$ is semisimple. 
\eeth

\bcor Let $V$ be a low height representation of $H$ and $v \in V$ an
element such that the $H$-orbit of $v$ in $V$ is closed. Then $V$ is a
semisimple representation for the reduced stabiliser $H_{v,red}$ of $v$.
\ecor

\bprop\label{cline} Let $H$ be as above and let $V$ be a low-height
representation of $H$. Then we have the following vanishing of group
cohomology:
\[
H^{i}(H,V) = 0
\]
for all $i \geq 1$.
\eprop
\noindent
{\it Proof.} We now recall from (\cite{mour} (pp 25,26) )the following
general result on low height modules of connected reductive groups:
Let $V$ be a low height module of $H$.  Let $\lambda$ be a dominant
weight which occurs in $V$. Then, if $V(\lambda) = H^{0}(\lambda)$ is
the dual of the Weyl module associated to $\lambda$, by the definition
of height and the low height property of $V$, it follows that
$V(\lambda)$ are also low-height $H$-modules.  In particular, it
follows that $V(\lambda)$ are also {\it irreducible} and they coincide
with their socle $L(\lambda)$. Therefore, by the semisimplicity of low
height modules, one has $V =~\stackrel{\lambda}\oplus V(\lambda)$.

Therefore by the Vanishing Theorem of Cline-Parshall-Scott-van der
Kallen (cf. \cite{jan} pp 237) we have the required cohomology
vanishing since
\[
H^{i}(H,V) =~\stackrel{\lambda}\bigoplus H^{i}(H,V(\lambda)) = 0
\]
for all $i \geq 1$. Q.E.D.

\subsection{\it Height and semistability}

 Let $F$ be a $G$-variety. Then a section $s : X \lr E(F)$ can be described
as a morphism from $\psi: E \lr F$ such that $\psi(e.g) = g^{-1}.s(e)$.  In
particular, if $H \subset G$ and $F = G/H$ then a section of $E(G/H)$ gives
a reduction of structure group of $E$ to $H$. 

We now recall the definitions of semistable, polystable and stable
principal bundles. Note that these definitions make sense for reductive
groups as well.

\bdefe (A. Ramanathan) $E$ is {\it semistable} (resp. {\it stable}) if for
every parabolic subgroup $P$ of $H$, and for every reduction of structure
group $\sigma_P : X \lr E(H/P)$ to $P$ and for any dominant character
$\chi$ of $P$, the bundle $\sigma_P^* (L_\chi))$ has degree $\leq 0$
(resp.$<0$).(cf.\cite {r1}).

\edefe

\bdefe\label{admis} A reduction of structure group of $E$ to a
parabolic subgroup $P$ is called {\it admissible} if for any character
$\chi$ on $P$ which is trivial on the center of $H$, the line bundle
associated to the reduced $P$-bundle $E_P$ has degree zero.
\edefe

\bdefe\label{poly} An $H$-bundle $E$ is said to be {\it polystable} if
it has a reduction of structure group to a Levi $R$ of a parabolic $P$
such that the reduced $R$-bundle $E_R$ is stable and the extended $P$
bundle $E_R(P)$ is an admissible reduction of structure group for $E$.
\edefe

\brem\label{hb} We note that there is a natural action of the group
$\Aut_GE$, of automorphisms of the principal $G$-bundle $E$, on $\Gamma (X,
E(G/H))$ and the orbits correspond to the $H$-reductions which are
isomorphic as principal $H$-bundles.

\erem

\brem\label{levired} Let $E_R$ be a stable $R$-bundle. Then  $E_R$ has no
further reduction of structure group to a Levi subgroup $L$ of a parabolic
subgroup in $R$.
\erem

\bprop\label{lh} Let $E$ be a principal $H$-bundle on $X$. Let $H
\hookrightarrow SL(V)$ be a {\it low height} faithful representation.
Then the following are equivalent: {\renewcommand{\labelenumi}{{\rm
(\alph{enumi})}}
\begin{enumerate}
\item The induced bundle $E(V) = E \times^H V$ is semistable.
\item The bundle $E$ is semistable as a principal $H$-bundle.

\end{enumerate}}
\eprop 
\noindent
{\it Proof.} (b) $\Rightarrow$ (a) follows by (\cite{imp} Theorem 3.1).

For (a) $\Rightarrow$ (b), we need to proceed as follows. By the Main
Theorem of \cite{mour}, any low height representation is actually
semi-simple.  Further, if $V = \bigoplus V_i$ is the decomposition
into irreducible $H$-modules, then the associated bundle $E(V)$
decomposes as $\bigoplus E(V_i)$ and the direct sum is of bundles of
degree zero.  Therefore it is clear that to prove the converse, we may
as well assume that the representation $\rho$ is an irreducible
representation of $H$ and also of low height. So since we are assuming
that $E$ is non-semistable, there exists a maximal parabolic subgroup
$P \subset H$ and a dominant character $\lambda$ such that the pull
back of $L_{\lambda}$ has degree $deg(L_{\lambda}) > 0$. Now it is not
very hard to see that there exists a parabolic $P_1$ in $SL(V)$ such
that $P = P_1 \cap H$ (cf \cite[Lemma 3.5]{imp}). Thus we see that,
the reduction of structure group of the vector bundle $E(V)$ to $P_1$
is given by a section $\sigma \in \Gamma(E(SL(V)/P_1))$ and the line
bundle $L_{\lambda}$ is the restriction of an ample line bundle
$L_{\lambda'}$ obtained by a dominant character $\lambda'$ of $P_1$.
It is clear that $deg(\sigma^*(L_{\lambda'})) > 0$ since
$degL_{\lambda} = deg(\sigma^*(L_{\lambda'}))$. This implies that
$E(V)$ is also non-semistable, and we are done.

\brem This theorem is strict in the sense that given an almost simple group
$H$ and a representation $H \lr SL(V)$ which is not of low height, there
exists a curve $X$ and a semistable $H$-bundle $E$ on $X$ such that $E(V)$
is not semistable. (The converse works for all but small primes. For more
precise details see \cite{imp}) \erem

\subsection{\it Functorial properness of the evaluation map}
The aim of this section is to define the basic functors and prove some
technical lemmas.
\noindent
Let $G$ be $SL(n)$ and let $H$ be a semisimple algebraic group, $H
\subset G$. For our convenience we make the following definition:

\bdefe\label{chevmod} Define an {\em affine $(G,H)$-module} $W$
associated to $(H,G)$ to be a finite dimensional $G$-module, such that
$G/H \stackrel{i}{\hookrightarrow} W$ is realised as a closed orbit of
a vector $w \in W$. Observe that since $G/H$ is affine, such a $W$
always exists. {\em We work with this canonical $W$ whenever we refer
  to the affine $(G,H)$-module}. This is a classical result ( cf. for
example \cite{dmil}, p 40 or \cite{bore}; also cf. Lemma \ref{chev}
below).  \edefe

\vspace{.2in}

Let
\[
F_G: (Schemes) \lr (Sets) 
\]
be the functor given by
\[
F_G(T) = \left \{ 
\begin{array}{l}
\mbox{isomorphism classes of semistable $G$-bundles of degree 0}\\
\mbox{on $X$ parametrised by $T$} \\
\end{array} \right.
\]
\noindent
One may similarly define the functor $F_H$ (note that since $H$ is
semisimple, for a principal $H$-bundle the associated vector bundles
have degree zero).

Let $x \in X$ be a marked point and let $F_{H,G,x}$ be the
functor
$$F_{H,G,x} (T) = \left\{
\begin{array}{l}
\mbox{isomorphism classes of pairs $(E,\sigma_x)$, $E = \{ E_t
\}_{t\in T}$} \\ 
\mbox{a family of semistable principal $G$-bundles of
degree 0}\\ 
\mbox{and $\sigma_x : T \lr E(G/H)_x$ a section}\\
\end{array} \right\}$$
(Recall that $E(G/H)_x$ denotes the restriction of $E(G/H)$ to $x \times
T \approx T$).

Notice that the functor $F_H$ is in fact realisable as the following
functor (by Remark \ref{hb}) .
$$F_{H,G}(T) = \left \{ \begin{array}{l} \mbox{isomorphism classes of
pairs $(E,s)$, $E = \{ E_t\}_{t \in T}$}\\ \mbox{a family of
semistable $G$-bundles of degree 0 and} \\ \mbox{$s = \{ s_t \}_{t \in
T}$ a section of E(G/H) on $X \times T$}\\ \mbox{or what we may call a
family of sections of $\{ E (G/H)_t \}_{t \in T}$}. \\
\end{array} \right \}$$
In what follows, we shall identify the functors $F_H$ with
$F_{H,G}$. With these definitions we have the following:

\bprop\label{dm}
Let $\alpha_x$ be the morphism induced by ``evaluation of
section'' at $x$:
\[
\alpha_x : F_H \lr F_{H,G,x}.
\] 
Then $\alpha_x$ is a {\it proper morphism of functors}.(cf. \cite{dm}).
\eprop
\noindent
{\it Proof.}  Let $T$ be an affine smooth curve and let $l \in T$. Let us
write $T^* = T-l$. Then by the valuation criterion for properness, we need
to show the following:

If $E$ is a family of semistable principal $G$-bundles on $X \times T$
together with a section $\sigma_x: T\lr E(G/H)_x$; such that for $t
\in T^*$, we are given a family of $H$-reductions, i.e. a family of
sections $s_{T^*} = \{ s_t \}_{t \in T^*}$, where, $s_t: X \lr
E(G/H)_t$, with the property that, at $x$, $s_t (x) = \sigma_x (t)$
$\forall ~ t \in T^*$; then we need to extend the family $s_{T^*}$ to
a section $s_T$ of $E(G/H)$ on $X \times T$ such that $s_l(x) =
\sigma_x (l)$ as well.

Let $W$ be an affine $(G,H)$-module associated to $(H,G)$.(see Def
\ref{chevmod}). Thus we get a closed embedding
$$E(G/H) \hra E(W)$$
and a family of vector bundles $\{ E(W)_t\}_{t\in
T}$ together with a family of sections $s_{T-l}$ and
evaluations $\{ \sigma_x (t) \}_{t \in T}$ such that $s_t(x)
= \sigma_x (t)$, $t \neq p$.

For the section $s_{T-l}$, viewed as a section of
$E(W)_{T-l}$ we have two possibilities:
\begin{enumerate}
\item[(a)] it extends as a regular section $s_T$.
\item[(b)] it has a pole along $X \times l$.
\end{enumerate}

Observe that if (a) holds, then we have
$$s_T (X \times (T-l)) \subset E(G/H) \subset E(W),$$ 
since
$E(G/H)$ is closed in $E(W)$, it follows that $s_T (X
\times l) \subset E(G/H)$.  Thus $s_l(X) \subset E(G/H)_l$.
Further by continuity, $s_l(x) = \sigma_x(l)$ as well and
this proves the proposition.

If (b) holds the reduction section exists over an open $U_T \subset X_T$
which contains all the primes of height 1 in $X_T$; or equivalently, the
$H$-bundle exists over $U_T$. We can now appeal to a theorem of
Colliot-Th\'el\`ene and Sansuc (\cite{coll} Theorem 6.13 pp 128) which
enables us to extend the principal bundle to $X_T$. In other words (b)
cannot occur. Finally observe that the limiting $H$-bundle is semistable
since it arises as a reduction of structure group of a semistable
$G$-bundle and $H \subset G$ is low height.  Q.E.D.

\brem A different proof involving the semistability of $E(W)$ is given in
\cite{bals}. Here we avoid it so as to improve the prime bounds arising out
of height considerations involved in the construction of the moduli spaces.
\erem

\section{Construction of the moduli space}

The aim of this section is to give a construction of the moduli space
of $H$-bundles. This section is somewhat different from the
corresponding one in\cite{bals} to enable us to provide the best prime
bounds.

Recall that $G = SL(n)$ and $H \subset G$ a semisimple subgroup.

We recall very briefly the Grothendieck Quot scheme used in
the construction of the moduli space of vector bundles
(cf. [Ses]).

Let $\cf$ be a coherent sheaf on $X$ and let $\cf (m)$ be
$\cf \otimes \co_X (m)$ (following the usual notations).
Choose an integer $m_0 = m_0 (n,d)$ ($n=$ rk, $d=$ deg) such
that for any $m \geq m_0$ and any semistable bundle $V$ of
rank $n$ and deg $d$ on $X$  we have
$h^i (V(m)) = 0$ and $V(m)$ is generated by its global
sections.

Let $\chi = h^0(V(m))$ and consider the Quot scheme $Q$
consisting of coherent sheaves $\cf$ on $X$ which are
quotients of $k^{\chi} \otimes_k \co_X$ with a fixed
Hilbert polynomial $P$.  The group $\cg = GL (\chi)$
canonically acts on $Q$ and hence on $X \times Q$ (trivial
action on $X$) and lifts to an action on the universal sheaf
$\ce$ on $X \times Q$.

Let $R$ denote the $\cg$-invariant open subset of $Q$
defined by 
\[
R = \left \{ q \in Q \mid 
\begin{array}{l}
\ce_q = \ce \mid_{X \times q} \mbox{
is locally free and the canonical map }\\
 k^{\chi}\lr H^0(\ce_q) \mbox{ is an isomorphism}, 
\det \ce_q \simeq \co_X
\end{array} \right \}
\]

We denote by $Q^{ss}$ the $\cg$-invariant open subset of $R$
consisting of semistable bundles and let $\ce$ continue to
denote the restriction of $\ce$ to $X \times Q^{ss}$.

Let $q'' :$ (Sch) $\lr$ (Sets) be the following functor:
$$q'' (T) = \left \{ (V_t,s_t) \mid \begin{array}{l}
\mbox{$\{ V_t\}$ is a family of semistable principal
$G$-bundles} \\
\mbox{parameterised by $T$ and $s_t \in \Gamma
(X,V(G/H)_t)~~\forall~ t \in T$}\\
\end{array} \right \}.$$
i.e. $q'' (T)$ consists pairs of rank $n$ vector bundles (or
equivalently principal $G$-bundles) together with a
reduction of structure group to $H$.

By appealing to the general theory of Hilbert schemes, one
can show without much difficulty (cf. [R1, Lemma~3.8.1]) that $q''$
is representable by a $Q^{ss}$-scheme, which we denote by $Q''$.

Let $W$ be an affine $(G,H)$-module associated to $(H,G)$ (see Def
\ref{chevmod}).  

\brem\label{J} Let $E \in Q^{ss}$ and consider $E(W)$ the associated vector
bundle.  Then, by the boundedness of the family $E \in Q^{ss}$ it follows
that there exists an $m_0$ ({\it independent of $E$}) such that if $s$ is
any section of $E(W)$, then ${\#}(zeroes (s)) \leq m_0$. Fix a subset $J
\subset X$ such that ${\#}J > m_0$.  \erem

The universal sheaf $\ce$ on $X \times Q^{ss}$ is in fact a
vector bundle. Let $\ce_G$ denote the associated
principal $G$-bundle, set 
\[
Q' = (\ce_G /H)_J = (\ce_G/H)_{x_j} \times_{Q^{ss}} (\ce_G/H)_{x_j'}
\cdot\cdot\cdot
\]
the fibre product being taken over all $j \in J$.
Then in our notation $Q' = \ce_G(G/H)_J$ i.e. we take the bundle over
$X \times Q^{ss}$ associated to $\ce_G$ with fibre $G/H$ and take its
restriction to $x_j \times Q^{ss} \approx Q^{ss}$ and take the product
over $Q^{ss}$.  Let $f: Q' \lr Q^{ss}$ be the natural map.  Then,
since $H$ is reductive, $f$ is an {\it affine morphism}.

Observe that $Q'$ parametrises semistable vector bundles
together with ``initial values of reductions to $H$ ''.

Define the ``evaluation map'' of $Q^{ss}$-schemes as follows:
$$\phi_J : Q''\lr Q'$$
$$(V,s) \longmapsto \{(V,s(x)) | x \in J \}.$$

\blem\label{propin} The evaluation map $\phi_J: Q'' \lr Q'$ is {\it
  proper} and {\it injective}.  \elem
\noindent
{\it Proof.}  Let $G/H \hookrightarrow W$ be as in Definition
\ref{chevmod} and let $(E,s)$ and $(E',s') \in Q''$ such that $\phi_J
(E,s) = \phi_J (E',s')$ in $Q'$.  i.e. $(E,s(x)) = (E',s'(x)) \forall
x \in J$.  So we may assume that $E \simeq E'$ and that $s$ and $s'$
are two different sections of $E(G/H)$ with $s(x) = s'(x) \forall x
\in J$.

Using $G/H \hookrightarrow W$, we may consider $s$ and $s'$ as
sections in $\Gamma (X,E(W))$ and further, as sections of $E(W)$ one
has $s(x) = s'(x) \forall x \in J$. By Remark \ref{J} this implies $s
= s'$. This proves that $\phi_J$ is {\it injective} (since $E(G/H)
\hookrightarrow E(W)$ is a closed embedding).

The properness of the map follows easily, the proof being as in
Proposition \ref{dm}. Thus $\phi_J$ being proper and injective is
{\it affine}. Q.E.D.

\brem In \cite{bals} a single base point served the purpose. Here we
employ the standard trick of increasing the number of base points to
achieve injectivity without the semistability of $E(W)$. This enables
us to improve the prime bounds.
\erem

\brem
It is immediate that the $\cg$-action on $Q^{ss}$ lifts to an action
on $Q''$.

Recall the commutative diagram
\begin{center}
\begin{picture}(400,325)(0,20)
\put(25,290){\ctext{$Q''$}}
\put(75,290){\vector(1,0){200}}
\put(325,290){\ctext{$Q'$}}
\put(325,250){\vector(0,-1){200}}
\put(75,265){\vector(1,-1){200}}
\put(325,25){\ctext{$Q^{ss}$}}
\put(125,150){\ctext{$\mu$}}
\put(350,150){\ctext{$f$}}
\put(250,200){\ctext{$\circlearrowright$}}
\put(175,320){\ctext{$\phi_J$}}
\end{picture}
\end{center}
By Lemma \ref{propin}, $\phi_J$ is a proper injection and hence affine.
One knows that $f$ is affine (with fibres $(G/H)^{\#J}$).  Hence $\mu$ is a
$\cg$-equivariant affine morphism.  \erem

\blem\label{cumrem} (cum remark) Let $(E,s)$ and $(E',s')$ be in the same
$\cg$-orbit of $Q''$.  Then we have $E \simeq E'$.  Identifying $E'$ with
$E$, we see that $s$ and $s'$ lie in the same orbit of $\Aut_GE$ on $\Gamma
(X,E(G/H))$.  Then using Remark \ref{hb}, we see that the reductions $s$
and $s'$ give isomorphic $H$-bundles.

Conversely, if $(E,s)$ and $(E',s')$ such that $E \simeq E'$ and the
reductions $s,s'$ give isomorphic $H$-bundles, using again Remark \ref{hb},
we see that $(E,s)$ and $(E',s')$ lie in the same $\cg$-orbit.  \elem

Consider the $\cg$-action on $Q''$ with the linearisation
induced by the {\it affine $\cg$-morphism} $\mu: Q'' \lr Q^{ss}$.  It is
seen without much difficulty that, since a good quotient of
$Q^{ss}$ by $\cg$ exists and since $Q''\lr Q^{ss}$ is an affine
$\cg$-equivariant map, a good quotient $Q''/\cg$ exists
(cf. [R1, Lemma~4.1]).

Moreover by the universal property of categorical quotients,
the canonical morphism
$$\ol{\mu} : Q''//\cg \lr Q^{ss}//\cg$$
is also {\it affine}.

Let $M_X(H)$ denote the scheme $Q''//\cg$.  then we have proved the
following theorem.  

\bth\label{moduli} Let $H \hookrightarrow G$ be a faithful low height
representation. Then there exists a coarse moduli scheme $M_X(H)$, for
semistable principal $H$-bundles. Further, the moduli space $M_X(H)$
is quasi-projective and the canonical morphism $\ol{\mu}: M_X(H) \lr
M_X(G)$ is {\it affine}.\eeth

\subsection{\it Points of the moduli}
In this subsection we will recall briefly the description of the
$k$-valued points of the moduli space $M_X(H)$. The general functorial
description of $M_X(H)$ as a coarse moduli scheme follows by the usual
process.

\bprop\label{points} The ``points '' of $M_X(H)$ are given by polystable
principal $H$-bundles.  \eprop

We firstly remark that since the quotient $q: Q'' \lr M_X(H)$ obtained
above is a good quotient, it follows that each fibre $q^{-1}(E)$ for
$E \in M_X(H)$ has a unique closed $\cg$-orbit. Let us denote an orbit
$\cg \cdot E$ by $O(E)$. The proposition will follow from the
following:

\blem If $O(E)$ is closed then $E$ is polystable.  
\elem 
\noindent
{\it Proof.}  Recall the definition of a polystable bundle Def
\ref{poly} and the definition of {\it admissible reductions} Def
\ref{admis}. If $E$ has no admissible reduction of structure group to
a parabolic subgroup then it is polystable and there is nothing to
prove. 

Suppose then that $E$ has an admissible reduction $E_P$, to $P \subset
H$.  Recall by the general theory of parabolic subgroups that there
exists a 1-PS $\xi :{\bf G}_m \lr H$ such that $P = P(\xi)$. Let
$L(\xi)$ and $U(\xi)$ be its canonical Levi subgroup and unipotent
subgroup respectively. The Levi subgroup will be the centraliser of
this 1-PS $\xi$ and one knows $P(\xi) = L(\xi) \cdot U(\xi) = U(\xi)
\cdot L(\xi)$. In particular, if $h \in P$ then $lim \xi(t) \cdot h
\cdot \xi(t)^{-1}$ exists. From these considerations one can show that there
is a morphism
\[
f : P(\xi) \times {\bf A}^{1} \lr P(\xi)
\]
such that $f(h,0) = m \cdot u$, where $h \in P$ and $h = m \cdot u$,
$m \in L$ and $u \in U$.  (see Lemma 3.5.12 \cite{r1})
 
Consider the $P$-bundle $E_P$. Then, using the natural projection $P
\lr L$ where $L = L(\xi)$, we get an $L$-bundle $E_P(L)$. Again, using
the inclusion $L \hookrightarrow P \hookrightarrow H$, we get a new
$H$-bundle $E_P(L)(H)$. Let us denote this $H$-bundle by $E_P(L,H)$.
It follows from the definition of admissible reductions and
polystability that $E_P(L,H)$ is {\em polystable}. 

Further, from the family of maps $f$ defined above, and composing with
the inclusion $P(\xi) \hookrightarrow H$ we obtain a family of
$H$-bundles $E_P(f_t)$ for $t \neq 0$ and all these bundle are
isomorphic to the given bundle $E$. Following (\cite{r1} Prop.3.5 pp
313), one can prove that the bundle $E_P(L,H)$ is the limit of
$E_P(f_t)$. It follows that $E_P(L,H)$ is in the $\cg$-orbit $O(E)$
because $O(E)$ is closed. Now by Lemma \ref{cumrem}, $E \simeq
E_P(L,H)$, implying that $E$ is polystable. Q.E.D.

\brem In the above Proposition we have only stated that there is a
surjective map from the set of isomorphism classes of polystable
$H$-bundles to the points of the moduli space. We believe that this
correspondence is a bijection but one possibly needs to discard a few
more primes.\erem

A few remarks are in order regarding the existence and properness of
the moduli spaces of principal bundles for ``large'' primes.

\brem\label{larhep} `` A general principle is that if a statement is true
in characteristic zero then it is also true for {\em large}
$p$''(cf. (\cite{mour})). One might therefore think that this would imply
the existence and projectivity of the moduli spaces of semistable
$H$-bundles for large primes, since one already knows this in char 0
(cf. for example \cite{r1}, \cite{f} or \cite{bals}).

We observe that this principle would indeed hold if one could show that the
subset corresponding to the semistable bundles in a family of $H$-bundles
is open over ({\bf Z} or large p); for the moduli spaces of $H$-semistable
bundles is realised as a GIT-quotient of a quasi-projective scheme and the
required results would follow by ``reduction modulo $p$'' for large $p$. To
the best of our knowledge the required ``openness result'' does not follow
by any general principle.

A key point of this paper is that even for the existence of the moduli
spaces of semistable $H$-bundles as a quasi-projective scheme for large
$p$, one requires {\it height} considerations. Moreover we give explicit
bounds for $p$.

Once the moduli space exists as a quasi-projective scheme for large $p$,
its projectivity follows for an unspecified {\it larger $p$}. One of the
hard parts of this paper is to give specific representation theoretic bounds
for $p$ for the projectivity of the moduli spaces.

\erem

\section{Separable index and slice theorem}

Let $T$ be a torus and $W$ be a finite dimensional $T$-module.
Further, let ${X}(T)$ be the free abelian group of characters of $T$
and $\cal S$ be the set of distinct characters that occur in $W$.

For every subset $S \subset \cal S$ we have the following map:
\[
\nu_S: {\bf Z}^{|S|} \lr {X}(T)
\]
given by $e_s \lr \chi_s$.

Let $g_S$ be the g.c.d of the maximal minors of the map $\nu_S$
written under the fixed basis. For any vector  $w\in W$, consider the
subset $S_w \subset \cal S$, consisting of characters that occur in $w$
with nonzero coefficients. i.e., if $w = \sum a_{\chi}(w)e_{\chi}$, then
\[
S_w = \{{\chi \in \cal S} | a_{\chi}(w) \neq 0 \}
\]
Then we have the following:

\blem The characteristic of the field, $p$ does not divide $g_{S_w}$ if 
and only if the action of $T$ on the vector $w$ is separable.  
\elem
\noindent
{\em Proof.} Let $T\cdot w$ denote the orbit of $w$ under the $T$
action.  Let $T_w$ be the stabilis
er.  Then $T/T_w$ is a torus and the
character group ${X}(T_w)$ of the stabiliser is the quotient of the
character groups ${X}(T)/{X}(T/T_w)$.  Moreover the image of the dual
map of the quotient map $T\to T/T_w$ is canonically identified with
the image of $\nu_{S_w}$.  Hence the group $T_w$ is reduced if and
only if this cokernel, identified with the cokernel of $\nu_{S_w}$,
has no $p$-torsion.  But this cokernel has $p$-torsion if and only if
the rank of $\nu_{S_w}$ drops mod $p$, which in turn happens if and
only if $p$ divides all the maximal minors. Hence the lemma.

\noindent
{\bf Notation}
\[
{\goth p}_T(W) = \{{\rm largest~prime~which~divides}~g_S |\forall S
\subset \cal S\} 
\]
\bdefe\label{sepind} Let $H \lr SL(W)$ be a finite dimensional
representation of $H$. Define the {\it separable index}, $\psi_H(W)$
of the representation as follows:
\[
\psi_H(W) = max\{ht_H(W), ~{\goth p}_T(W)\}
\]
\edefe

\brem When $T$ is a maximal torus of a semisimple group $H$ and the
$T$-module $W$ is actually an $H$-module, then the set of characters
that occur on $W$ can be written down explicitly using Standard
Monomial Theory. From this very explicit form, this separable index is
computable, though it could be tedious or may need a computer. The few
cases where we made some computations indicated that this index is
possibly bounded above by the dimension of $W$. One can easily observe
that the absolute value of each minor of the map $\nu_S$ is bounded
above by $l! \cdot h^{l}$, where $l = rank(G)$ and $h =
ht_G(W)$. Hence the separable index has a weak upper bound given by
$l!  \cdot h^{l}$.\erem

\bdefe A representation $H \lr SL(W)$ is said to be with {\bf low
  separable index} if $p > \psi_H(W)$.  \edefe

\bth\label{strongsep} If $W$ is a low separable index $H$-module then the
action of $H$ on $W$ is {\it strongly separable} i.e., the stabilizer at
any point is absolutely reduced.  
\eeth 
\noindent
{\em Proof.} Since the representation is low height, every nilpotent
in the Lie algebra of the $H$ is integrated in $SL(W)$ and hence the
nilpotent part of the Lie algebra of the stabiliser at any $w \in W$
will actually lie in the Lie algebra of the reduced stabiliser. Thus
by the Jordan decomposition of the stabiliser, it is enough to ensure
separability of the action of any maximal torus. Separability index
assures that the given maximal torus $T$ acts separably at all points in
$W$. This implies that every maximal torus acts separably at all points as
all maximal tori are conjugates. Hence the action of $H$ is strongly
separable.  Q.E.D.

\brem\label{saturate} We recall briefly the notions of {\it saturated
subgroups} of $GL(V)$. For details cf.pp 524-526 \cite{tens}.

We first define a one parameter subgroup defined by an element of
order $p$. Let $V$ be a finite dimensional $k$-vector space, and let
$s \in GL(V)$ be an element such that $s^{p} = 1$. One has $s = 1 + u$
where $u^p = 0$. If $t \in k$, we can define an element $s^{t} \in
GL(V)$ by the truncated binomial formula:
\[
s^{t} = 1 + tu + \frac {t(t-1)u^2}2 + ...
\]
summed upto $u^i$ with $i < p$.
The map $t \lr s^t$ defines a homomorphism of algebraic groups:
\[
\phi_s :{\bf G}_a \lr GL(V)
\]
where ${\bf G}_a$ is the additive group. This homomorphism has two characterising properties:
\begin{itemize}
\item $\phi_s(1) = s$
\item The map $t \lr \phi_s(t)$ is a polynomial of degree $< p$.
\end{itemize}
Let $H\subset GL(V)$ be a subgroup. We say that $H$ is {\it saturated} if every unipotent element $s \in H$ has the following properties: 
\begin{itemize}
\item $s^p = 1$
\item $s^t \in H$ for every $t \in k$
\end{itemize}
One can see that given any subgroup $H$ there is a smallest saturated
subgroup which contains $H$ called the {\it saturation of $H$}.

A property of saturated groups which we need is that if $H$ is
saturated and $H^0$ is the connected component of identity of $H$ then
the index $[H:H^0]$ is coprime to $p$. (cf.pp 524-526 \cite{tens}).

One can again generalize all these notions for an arbitrary reductive
algebraic group $G$ instead of $GL(V)$. Among the elementary examples
of saturated subgroups are parabolic subgroups, centralizers of any
subgroup, and Levi subgroups (since they can be realised as the
centralizer of a torus). We can isolate a couple of key properties in
the theory of low height representations:
\begin{enumerate}
\item If $G \lr GL(V)$ is a
low height representation of $G$ then the isotropy subgroups of closed
orbits in $V$ are saturated. 
\item If $S$ is a reductive and saturated subgroup of $G$ and if $G
\lr GL(V)$ is a low height representation of $G$ then $V$ is a low
height module for $S$ as well.(cf. \cite{mour} p.25)
\end{enumerate}
\erem

\bprop\label{luna}(A version of Luna's \'etale slice theorem in
char.$p$) Let $W$ be a low separable index $G$-module. Let $F$ be a
fibre of the good quotient $q :W \lr W//G$, and let $F^{cl}$ be the
unique closed orbit contained in $F$. Then there exists a $G$-map
\[
F \lr F^{cl}.
\]
\eprop
\noindent
{\em Proof.} Since $\psi_G(W) = max\{ht_G(W), ~{\goth p}_T(W)\}$ the
assumption $p > \psi_G(W)$ on the separable index implies the
following:
\begin{enumerate}
\item Every stabiliser subgroup for the $G$-action on $W$ is reduced,
the action being {\it strongly separable} (by Th \ref{strongsep}).
\item It is {\it saturated}, the representation $G \lr GL(W)$ being
low height.
\item When $w$ is a quasistable point in $W$ or equivalently, the orbit
$G \cdot w$ is closed, then $W$ is a semisimple representation of the
stabiliser $G_w$. This is a consequence of the main theorem of
\cite{tens}, namely that low height representations are semi-simple.
\end{enumerate}

For more on this (cf. \cite{mour} pp 20-25); it may be kept in mind
that the height of the representation, $ht_G(W)$, coincides with
Serre's index $n_G(W)$.

A close examination of Luna's proof shows that the key point is the
complete reducibility of the tangent space $T_w(W)$, of the affine
$G$-module. This is used then to get a splitting of the canonical
injection of the tangent space of the closed $G$-orbit in
$T_w(W)$. Once this is achieved the slice can be constructed. The
above proposition is then a corollary to the main slice theorem
applied to a single orbit. (For details cf. \cite{bard} Prop 8.5 p
312). Q.E.D.

\section{Towards the flat closure}
Fix as in \S3.2 a faithful {\it low height} representation $H
\hookrightarrow G$ defined over $k$ as well as an {\em affine
$(G,H)$-module} associated to the pair $(H,G)$.(cf.Def \ref{chevmod}).

Consider the extension of structure group of the bundle $P_K$ via the
induced $K$-inclusion $H_K \hookrightarrow G_K$ . We denote the
associated $G_K$-bundle $P_K(G)$ by $E_K$.

Then, since $G = SL(n)$, by the projectivity of the moduli space of
 semistable vector bundles, there exists a {\it semistable extension}
 of $P_K(G) = E_K$ to a $G_A$-bundle on $X \times \spec A$, which we
 denote by $E_A$.  Call the restriction of $E_A$ to $X \times l$
 (identified with $X$) the {\it limiting bundle} of $E_A$ and denote
 it by $E_l$ (as in \S1). One has in fact slightly more, which is what
 we need.

\blem\label{langt} Let $E_K$ denote a family of semistable $G_K$-bundles on
$X \times \spec K$ (or equivalently a family of semistable vector bundles
of rank $n$ and trivial determinant on $X \times T^*$). Then (by going to
a finite cover $S$ of $T$ if need be ) the principal bundle $E_K$ extends
to $E_A$ with the property that the limiting bundle $E_l$ is in fact {\it
polystable} i.e, a direct sum of stable bundles.  \elem
\noindent
{\it Proof.}  The proof of this Lemma is possibly well known but for
the sake of completeness we give it here. Recall notations as in \S4
regarding Quot schemes etc.

Note that the moduli space in question, namely of semistable principal
$G$-bundles, is a GIT quotient $Q^{ss} \longrightarrow M$ by ${\cal
G}$, and the family $E_A(G)$ is given by a morphism $T \longrightarrow
M$.  Lift the $K$-valued point, namely, $r_K$, given by the family
$E_K$, to $Q^{ss}$ and consider the ${\cal G}$-orbit $R_0$ of $r_K$ in
$Q^{ss}$. Let $\overline R_0$ be its closure in $Q^{ss}$. Since the
$K$-valued point $r_K$ is in fact an $A$-valued point of $M$, the GIT
quotient of $\overline R_0$ is indeed the curve $T$. Also, observe
that the closure intersects the closed fibre. Consider the morphism
$\psi : \overline R_0 \lr T$. Since the base is a curve $T$, one has a
{\it multi-section} for the morphism $\psi$, and one obtains the curve
$S$. The general fibre has been modified only in the orbit, therefore
the isomorphism class of the bundles remains unchanged. Q.E.D.

\brem It is to be noted that the definition of polystability given
here coincides with that in Def \ref{poly}, in the sense that a closed
orbit in the Quot scheme corresponds to a polystable vector bundle.
\erem

We observe the following:
\begin{itemize}
\item Note that giving the $H_K$-bundle $P_K$ is giving a reduction of
 structure group of the $G_K$-bundle $E_K$ which is equivalent to
 giving a section $s_K$ of $E_K (G_K/H_K)$ over $X_K$ .

\item We fix a base point $x \in X$ and denote by $x_A = x \times
\spec A$, the induced section of the family (which we call the {\it
base section}):
\[
X_A \lr \spec A
\]

\item Let $E_{x,A}$ (resp $E_{x,K}$) be as in \S1, the restriction of
$E_A$ to $x_A$ (resp $x_K$). Thus, $s_K(x)$ is a section of
$E_K(G_K/H_K)_x$ which we denote by $E_x (G_K/H_K)$.

\item Since $E_{x,A}$ is a principal $G$-bundle on $\spec A$ and
therefore trivial, it can be identified with the group scheme $G_A$
itself. {\it For the rest of the article we fix one such
identification, namely:}
\[
\xi_A : E_{x,A} \lr G_A.
\]
\item Since we have fixed $\xi_A$ we have a canonical identification
\[
E_x(G_K/H_K) \simeq G_K/H_K 
\]
which therefore carries a natural {\it identity section} $e_K$ (i.e
the coset $id.H_K$). Using this identification we can view $s_K(x)$ as
an element in the homogeneous space $G_K/H_K$.

\item Let $\theta_K \in G(K)$ be such that $\theta_K^{-1} \cdot s_K(x)
= e_K$.  Then we observe that, the isotropy subgroup scheme in $G_K$
of the section $s_K(x)$ is $\theta_K.H_K.\theta_K^{-1}$. 

\item On the other hand one can realise $s_K(x)$ as the identity coset of
$\theta_K.H_K.\theta_K^{-1}$ by using the following
identification:
 $$G_K/\theta_K.H_K.\theta_K^{-1} 
 \stackrel{\sim}{\lr} G_K/H_K.$$
 $$g_K(\theta_K.H_K.\theta_K^{-1}) \longmapsto
 g_K\theta_K.H_K.$$ 
\end{itemize}
\noindent 
\bdefe Let ${H_K'}$ be the subgroup scheme of $G_K$ defined as: 
\[
H_K' := \theta_K.H_K.\theta_K^{-1}.
\]
\edefe
\noindent
Using $\xi_A$ we can have a canonical identification:
\[
 E_{x}(G_K/H_K') \simeq G_K/H_K'.
\] 
\noindent
Then we observe that, using the above identification we get a section
$s_K'$ of $E_K (G_K/H_K')$, with the property that, $s_K'(x)$ is the
{\it identity section} and moreover, since we have conjugated by an
element $\theta_K \in G_A(K) (=G(K))$, the isomorphism class of the
$H_K$-bundle $P_K$ given by $s_K$ does not change by going to $s_K'$.

Thus, in conclusion, the $G_A$-bundle $E_A$ has a reduction to
$H_K'$ given by a section $s_K'$ of $E_K (G_K/H_K')$, with the
property that, at the given base section $x_A = x \times \spec A$, we
have an equality $s_K' (x_A) = e_K'$, with the {\it identity element}
of $G_K/H_K'$ (namely the coset $id.H_K'$).

\bdefe\label{fc} The {\it flat closure} of the reduced group scheme
$H_K'$ in $G_A$ is defined to be the schematic closure of $H_K'$ in
$G_A$ with the reduced scheme structure.  Let $H_A'$ denote the
flat-closure of $H_K'$ in $G_A$.(cf. Lemma \ref{chev})
 \edefe
\noindent
We then have a canonical identification via $\xi_A$:
\[
 E_{x}(G_A/H_A') \simeq G_A/H_A'.
\] 
\noindent
One can easily check that $H_A'$ is indeed a subgroup scheme of $G_A$
since it contains the identity section of $G_A$, and moreover, it is
faithfully flat over $A$. Notice however that $H_A'$ {\it need not} be
a {\it reductive} group scheme; that is, the special fibre $H_k$ over
the closed point need not be reductive.

Observe further that $s_K'(x)$ extends in a trivial fashion to a
 section $s_A'(x)$, namely the {\it identity coset section} $e_A'$ of
 $E_x(G_A/H_A')$ identified with $G_A/H_A'$ .
 
 \brem\label{bull} If $H_A'$ is {\it reductive} then the semistable
 reduction theorem (Theorem \ref{ssred}) follows quite easily. Indeed,
 firstly by the {\it rigidity} of reductive group schemes over $\spec
 A$ (SGA 3, Expose III Cor 2.6 pp 117), by going to a finite cover, we
 may assume that $H_A' = H \times \spec A$. Secondly, in this case one
 can realise $H_A'$, as the isotropy subgroup scheme for a closed
 orbit $w_A \in W_A$.  Then we have a {\it closed $G$-immersion} of
 $G/H$ in a $G$-module $W$, and one may view $s_K$ as a section of
 $E_K(W_K)$.  Note that $ E_K (G_K/H_K') \subset E_K (W_K)$.

By choice, along $x_A$, the section $s_K (x)$ extends regularly to a
 section of $E_A (G_A/H_A') \subset E_A (W_A)$.  Hence by Proposition
 \ref{dm}, $s_K$ extends to a section $s_A$ which gives the required
 reduction over $X \times \spec A$.  \erem

\section{Affine embedding of $G_A/H_A'$}

As we have noted, $H_A'$ need not be reductive and the rest of the proof is
to get around this difficulty.  Our first aim is to prove that the
structure group of the bundle $E_A(G_A)$ can be reduced to $H_A'$ which is
the statement of Proposition \ref{flatext}.

 We need to prove the following generalisation of a well-known result
 (cf. for example \cite{bore}).

\blem\label{chev} 
There exists a finite dimensional $G_A$-module $W_A$ such that
 $G_A/H_A'~\hookrightarrow~W_A$ is a $G_A$-immersion.  
\elem
\noindent
{\it Proof.}  We follow the standard proof.  Let $I_K$ be the ideal
 defining the subgroup scheme $H_K'$ in $K(G)$ (note that $G_A$ (resp
 $G_K$) is an affine group scheme and we denote by $A(G)$ (resp
 $K(G)$) its coordinate ring).

Set $I_A = I_K \cap A(G)$. Then it is easy to see that since we are over a
discrete valuation ring, $I_A$ is in fact the ideal in $A(G)$ defining the
flat closure $H_A'$. Observe also that $I_A$ is a {\it primitive} $A$
submodule of $A(G)$, that is, $A(G)/I_A$ is torsion free; further, $I_A
\otimes k = I_k$ is the defining ideal in $k(G)$ of $H_k'$ in $G_k$ and
$I_A \otimes K$ is $I_K$. 

Since $A(G)$ and the other modules involved are free over the discrete
valuation ring $A$, a set generates $I_A \otimes k = I_k$ if and only
if it generates $I_A$. Thus we may now choose a finite generating set
$\{f_{i}\}$ of $I_A$, such that their images $f_{i,k}$ generate $I_k$.

As in the classical proof, one has a finite dimensional $G_K$-submodule,
$V_K$, containing the $\{f_{i}\}$. Now set $V_A = V_K \cap A(G)$ and $M =
V_A \cap I_A$. Observe that $I_A$, $V_A$ and hence $M$ are all
$G_A$-submodules of $A(G)$. This can be seen by keeping track of the
co-module operations. Then clearly $V_A$ is primitive in $A(G)$ and $M$ is
also primitive in $A(G)$ and in particular, primitive in $V_A$. 

If we set
\[       
M_k = M \otimes k~ \mbox{and}~V_k = V_A \otimes k
\]
 we see that the inclusion $M \hookrightarrow V_A$ induces an inclusion
 $M_{k} \hookrightarrow V_k$ . Observe that 
\[
f_{i}~\in M~, f_{i,k}~\in~ M_{k}~ \mbox{and}~ M~\subset~I_A
\]
\[
M_{k}~\subset~I_k~ \mbox{and}~ M_{k}~=~V_k~\cap~I_k
\]
 We claim that, for $g \in G_A(k)$, one has
\[
g \cdot M_{k} \subset M_{k} \iff g \in H_k'
\]
 Obviously, if $g \in H_k'$, then $g \cdot M_{k} \subset M_{k}$,
 since $V_k$ is $G$-stable and $I_k$ is $H_k'$-stable. Thus, it
 suffices to show that
\[
f_{i,k}(g) = 0~\mbox{for all}~ i 
\]
 that is,
\[
{f_{i,k}}\mbox{ vanish on }~g
\]
 Since $f_{i,k} \in M_{k}$, it suffices to show that 
\[
\phi(g) = 0~ \mbox{for}~ \phi~\in~M_{k}
\]
 But $\phi(g) = (g^{-1} \cdot \phi)(id)$, where $g^{-1} \cdot \phi$ is the
 action of $G$ on functions on $G$.  Now, by hypothesis, $(g^{-1} \cdot
 \phi) \in M_{k}$. Since $M_{k} \subset I_k$, and $id \in H_k'$ , we
 see that $(g^{-1} \cdot \phi)(id) = 0$ . This proves the above claim.

 Similarly, if we set 
\[
M_{F} = M \otimes_{A} F~ \mbox{and}~ V_F = V_A \otimes_{A} F
\]
 where $F$ is any field containing $A$, we see that for $g \in G(F)$
\[
g \cdot M_{F} \subset M_{F} \iff g \in H_A'(F)
\]
 Let $L$ denote the primitive rank one $A$-submodule $\wedge^d M
 \hookrightarrow \wedge^d V = W_A$, and [$L$] the $A$-valued point of
 ${\bf P}(W_A)$ defined by $L$. Here, ${\bf P}(W_A)$ is defined by the
 functor associated to rank one direct summands of $W_A$.  Then, the
 above discussion means that, we can recover $H_A'$ as the isotropy
 subgroup scheme at [$L$] for the $G_A$-action on ${\bf P}(W_A)$.

 Recall that, for any field $F$, the isotropy subgroup of $G_A(F)$, at
 the point of ${\bf P}(W_A (F))$ represented by the base change of $L$ by
 $F$, is $H_A'(F)$ .
 
 Fix a generator $l \in L$ so that $l$ is a primitive element in $W_A$
 and consider the isotropy subgroup scheme $H_A''$ at $l$ for the
 $G_A$-action on $W_A$.  We claim that, $H_A''$ coincides with $H_A'$.
 To see this, observe that, $H_A''$ is the subgroup scheme of $G_A$
 which leaves the closed subscheme (identified with $Spec(A)$)
 determined by $l$ invariant (with the corresponding automorphism on
 this subscheme being identity). We see then that, $H_A''$ is a {\it
   closed} subgroup scheme of $G_A$. Further, we see that since
 $H_A''$ is the isotropy subgroup of the vector $l\in L$ and $H_A'$
 that of the line [$L$] we have $H_A'' \hookrightarrow H_A'$.  Since
 $H_K'$ is semi-simple, it has no characters and therefore, the
 isotropy subgroup scheme at $(l \otimes K) \in (W_A \otimes K)$ is
 precisely $H_K'$. This means that, $H_K'' = H_K'$.  Now, $H_K'$ is
 open (dense) in $H_A'$ (since $H_A'$ is the flat closure of $H_K'$ in
 $G_A$ ) so that, $H_K'$ is also dense in $H_A''$.  This implies that,
 $H_A'$ and $H_A''$ coincide set theoretically.  Observe also that
 $H_A'$ is {\it reduced} by the definition of flat closure. Thus, it
 follows that $H_A'$ = $H_A''$. This implies that, $G_A/H_A'
 \hookrightarrow W_A$ is a $G_A$-{\it immersion} and the above lemma
 follows. Q.E.D.
 
 \brem Regarding the Lemma \ref{chev} proved above, we note that
 usually the subgroup scheme $H_A'$ can be realised only as the
 isotropy subgroup scheme of a line in a $G_A$-module. But here, since
 the generic fibre of $H_A'$ is semisimple, one is able to realise
 $H_A'$ as the isotropy subgroup scheme of a primitive element in a
 $G_A$-module and the limiting group scheme also as an isotropy
 subgroup scheme for an element in a $G_k$-module. We note here that
 last part of the above proof is seen easily by observing that a
 non-trivial character of $H_A'$ by definition is a non-trivial
 character of $H_K'$ and hence $H_A'' = H_A'$. \erem

 \brem\label{chev1} We make the following key observations about the
 group scheme $H_A'$. The flat group scheme $H_A' = Stab(w_A)$, is the
 isotropy subgroup scheme of $G_A$ at an $A$-valued point $w_A \in
 W_A$, where $W_A$ can be realised as $W \otimes A$ (after going to a
 finite cover of $A$ if need be) and $W$ is the affine $(G,H)$-module
 such that $G/H \subset W$.
 
 Moreover, it is also shown as a part of the proof that the closed
 fibre $H_k' = Stab(w_k)$, is the isotropy subgroup scheme of $G_k$
 for a vector $w_k \in W$.
 
 Thus if we assume that $p > \psi_G(W)$, it follows by Theorem
 \ref{strongsep} that $H_k'$ is {\bf reduced}.  \erem

\section{Semistable bundles, semistable sections and saturated groups}
The aim of this section is to prove some general lemmas on polystable
bundles and semistable sections. We assume that $p > \psi_G(W)$,
notations as in \S5.

\bdefe\label{bogo} (following Bogomolov) Let $E$ be a principal
$G$-bundle and let $G \lr GL(V)$ be a representation of $G$. Let $s$
be a section of the associated bundle $E(V)$. Then we call the section
$s$ stable (resp semistable, unstable) relative to $G$ if at one point
$x \in X$ (and hence at every point on $X$) the value of the section
$s(x)$ is stable (resp semistable, unstable).  \edefe

\noindent
(It is easy to see the non-dependence of the definition on the point $x
\in X$. Consider the inclusion $k[V]^{G} \hookrightarrow k[V]$ and the
induced morphism $V \lr V/G$. This induces a morphism $E(V) \lr
E(V/G)$. Observe that $V/G$ is a trivial $G$-module. Thus we have the
following diagram:
\[
s:X \lr E(V) \lr E(V/G) \simeq X \times V/G
\]
Composing with the second projection we get a morphism $X \lr V/G$
which is constant by the projectivity of $X$. Hence the value of the
section is determined by one point in its $G$-orbit.)  (cf. \cite{rou}
1.10)

\blem\label{stab} Let $E(W)$ be a semistable vector bundle of {\em
  degree zero} and let $R$ be a {\em saturated reductive} subgroup of
$GL(W)$. Suppose that $E(W)$ has a reduction of structure group to
$E_R$, a stable $R$-bundle, and further suppose that we have a
non-zero section $s: X \lr E_R(W) = E(W)$.  Then $s$ is a semistable
section in the sense of Def \ref{bogo}.  \elem
\noindent
{\it Proof.} Suppose that this is not the case. Then as observed in
the definition, if $s(x)$ unstable for a single $x \in X$ implies
it is unstable for all $x \in X$. In particular for the generic point
$x_0 \in X$. (cf. \cite{rou} Prop 1.5)

Since $s$ is a non-zero section of $E_R(W)$ and $E_R(W)$ is semistable
of degree zero, it is nowhere zero. This section gives a reduction of
structure group of $E_R(W)$ to a maximal parabolic subgroup $P_s$,
given by the extension:
\[
0 \lr {\cal O}_X \lr E_R(W) \lr V \lr 0
\]
for some degree zero vector bundle $V$ and where the first inclusion
is given by the section $s$.

Notice that $SL(W)/P_s = {\bf P}(W)$. In the language of \cite{rr}, the
section $s$ can be thought of as taking values in the cone $W$
and $deg (s) = 0$.

{\it We now claim that w.l.o.g we may assume that the representation
  $W$ is an irreducible $R$-module}.

Since $W$ is a low height $R$-module it is completely reducible, i.e
it can be expressed as
\[
W = \bigoplus W^{\alpha}
\]
where $W^{\alpha}$ are irreducible $R$-module. Any element $w \in W$ can be
expressed as $w = \oplus w^{\alpha}$ with $w^{\alpha} \in W^{\alpha}$. It is easy to see
that if, $w$ is $R$-unstable and if $\lambda$ is a Kempf 1-PS in $R$
which drives $w$ to 0 then $\lambda$ drives all the $w^{\alpha}$'s to 0 as
well. Further, as bundles
\[
E_R(W) = \bigoplus E_R(W^{\alpha})
\]
and since $E_R(W)$ is semistable of degree 0 all the $E_R(W^{\alpha})$
are semistable of degree 0 being direct summands of $E_R(W)$. The
given section also breaks up as $s = \oplus s^{\alpha}$ to give
non-zero (and hence nowhere zero!)  sections of $E_R(W^{\alpha})$
(since $s(x) = w = \oplus w^{\alpha}$, here of course, not all
$\alpha$'s may be involved!).  

Again by Def \ref{bogo}, the new sections $s^{\alpha}$ continue to
remain unstable since instability is determined at a point $x \in X$.
This proves the claim.

Once $W$ is irreducible as an $R$-module by Schur Lemma the connected
component $Z^0(R)$ of center of $R$ acts as scalars on $W$ and hence
{\it trivially} on ${\bf P}(W)$ and as scalars on the ample line
bundle $L$ on it. 

Since $m = s(x_0)$ is unstable we have a Kempf instability flag $P(m)$
and the corresponding 1-PS $\mu$, are also defined over the field
$K(X)$ . This follows by the low separable index assumption, namely $p
> \psi_G(W)$, which in particular implies $W$ is a low height module
for $G$ and hence for the saturated subgroup $R$ (cf. \cite[Prop
3.13]{rr} and \cite[Theorem 3.1]{imp}).

The parabolic subgroup being defined over $K(X)$ gives a reduction of
structure group of $E_R$ to a parabolic $P$ of $R$. Let $W = \bigoplus
W_i$ be the weight space decomposition of $W$ with respect to $\mu$.
Let $m = m_0 + m_1$, with $m_0$ of weight $j > 0$ and $m_1$ the sum of
terms of higher weights. In other words, in the projective space ${\bf
P}(W)$ we see that $\mu(t) \cdot m \lr m_0$. It is not too hard to see
that we have an identification of the Kempf parabolic subgroups
associated to the points $m$ and $m_0$. i.e $P(m) =
P(m_0)$.(cf. \cite[Proposition 1.9]{rr}).

In the generic fibre $E_R(W)_{x_0}$ we have the projection
\[
\bigoplus_{i \geq j} W_i \lr W_j
\]
which takes $m$ to $m_0$. This gives a line sub-bundle $L_0$ of degree
zero of $E_R(W)$ corresponding to $m_0$. It then follows that $m_0$ is
in fact semistable for the action of $P/U$, the Levi of $P$, for a
suitable choice of linearisation obtained by twisting the action by a
dominant character $\chi$ of $P$. (This is essentially the content of
\cite[Prop.1.12]{rr} and we can apply it since we work in the degree
0.)

The semistability of the point $m_0$ with this new linearisation the
forces the degree inequality:
\[
deg(L_0 \otimes L(\chi)^{-1}) \leq 0
\]
But since $deg(L_0) = 0$, this implies $deg(L(\chi)) \geq 0$. This
contradicts the stability of $E_R$. Q.E.D.

\brem We note that the condition of semistability of the vector bundle
$E_R(W)$ is assumed here since \cite{imp} proves it only for
semisimple groups. But in the situation in which we need (cf. Prop
\ref{lc}) this condition automatically holds since we have the
following inclusion
\[
R \hookrightarrow G \hookrightarrow GL(W)
\]
and therefore $E_R(W) = E(W)$ and $E(W)$ is semistable since $W$ is a
low height representation of $G$.  \erem

\blem\label{redstab} Let $E_R$ be a stable $R$-bundle as above and let
$I$ be a saturated reductive (possibly non-connected) subgroup of $R$
such that $E_R$ has a reduction of structure group to $I$. Then the
reduced $I$-bundle is also stable.  \elem
\noindent
{\it Proof.} We first {\it claim} that $I$ is {\it irreducible} in
$R$: if not, then by the low height property there exists a parabolic
$P$ and a Levi $L$ in it such that $I \subset L$ and this is
irreducible. This gives a reduction of structure group of $E_R$ to $L$
and this again contradicts the stability of $E_R$, by Remark
\ref{levired}.

Now to prove the Lemma, suppose that the reduced $I$-bundle $E_I$ is
{\it not stable}. Then, $E_I$ has an reduction of structure group
$\sigma$, to a maximal parabolic $P \subset I$.  Observe that any
parabolic subgroup of a reductive algebraic group looks like
$P(\lambda)$ for a 1-PS $\lambda :{\bf G}_m \lr I$. Now consider
$P_R(\lambda)$ the induced parabolic in $R$.  Then, it is clear that
$P_R(\lambda)$ gives a reduction of structure group for $E_R$.

Notice that $P_R(\lambda)$ in $R$ may not be a maximal parabolic, but
there exists a maximal parabolic $Q$ containing it. Now note that
by the irreducibility of $I \subset R$ seen above, $Q \cap I$ is a
proper parabolic in $I$ and contains $P_I(\lambda)$. Therefore by the
maximality of $P_I$ it follows that $Q \cap I = P_I$.

Let $\chi$ be a dominant character of $P(\lambda)$ and let the induced
line bundle be $L_\chi$ such that $deg(\sigma^{*}(L_\chi)) \geq 0$.
Then, since $Q$ is a maximal parabolic a multiple of $\chi$ extends to
a dominant character of $Q$ and the induced line bundle $L_\chi$ on
$I/P$ is the restriction of the line bundle from $R/Q$. Therefore, the
degrees of the pull backs to $X$ remain the same. This contradicts the
stability of $E_R$. Q.E.D.

\bprop\label{monodr} Let $E_R$ be a stable $R$-bundle and $s$ be a
non-zero section of $E_R(W)$ as in Lemma \ref{stab}. Let $s(x) = w$.
Then the $R$-orbit of $w$ is closed and $s$ takes its image in the
closed orbit.  \eprop
\noindent
{\it Proof.} By Lemma \ref{stab}, since $E_R$ is stable, $w \in
W^{ss}$. Therefore the section $s$ which can be thought of as a map
\[
s: E_R \lr W^{ss}
\]
which further takes its values in a fibre $F$ of the GIT quotient:
\[
W^{ss} \lr W^{ss}//R
\]
Thus the section $s$ gives the following map:
\[
s: E_R \lr F
\]
and $F$ contains the vector $w$. 

\noindent
{\it We need to show that the orbit $R \cdot w$ is closed}.
\noindent
We prove this by contradiction.

Suppose then that orbit of $R \cdot w$, is {\it not} closed. Let $I$ be the
isotropy at a point $f \in F$ such that $R \cdot f$ is closed.  Note that
the identity component $I^o$ is reductive and saturated and $I$ is
also reduced.  

Then by Proposition \ref{luna} we have an $I$-invariant ``slice'', $S \subset
F$ and an $R$-isomorphism
\[
\theta :R \times^{I} S \simeq F
\]
\[
\theta([r,s]) = r \cdot s
\]

This gives a $R$-equivariant morphism
\[
l: F \simeq R \times^{I} S \lr R/I \simeq F^{cl}.
\]
The composition $l \circ s = s_1$ of the maps $s$ and $l$ gives a
reduction of structure group, $E_I \subset E_R$ to the isotropy $I =
Stab_R(f)$ of a point $f \in F^{cl}$. 
\noindent
By Lemma \ref{stab} the $I$-bundle $E_I$ is stable.

Consider the given section $s$ of $E_R(W)$ as obtained via the
reduction of structure group to $I$. This is given as follows:
\[
s_1 :E_I \lr F \hookrightarrow W
\]
which is $I$-equivariant. Observe that without loss of generality (by
taking a conjugate of the isotropy $I$) we may assume that $w \in
Im(s_1)$.

(This is easy to see. Indeed, starting with a pair $(I,S)$ namely a
slice and an isotropy subgroup at $f \in S$, the given point $w \in F$
can be expressed as an equivalence class $ w = [r,s_0]$. Then by
translating the slice $S$ by the element $r \in R$ we get a new slice
$r \cdot S = S'$ and a new pair $(I',S')$ where $I' = r \cdot I \cdot r^{-1}$. It is
clear that we have an isomorphism
\[
F \simeq  R \times^{I} S \simeq R \times^{I'} S'
\]
and under this identification we get a reduction of structure group to
$I'$ with the property that the image of the section contains
the given vector $w$.)

Further, by assumption $w \in W - W^{I}$.

Moreover, the $I$-orbit closure of $w$ contains $f \in W^{I}$.
Therefore, if $\overline w$ is the image of $w$ in the quotient space
$W/W^{I}$, then clearly $\overline w$ is an $I$-{\it unstable} vector
in $W/W^{I}$.

Observe also that since $I$ is saturated, by \cite{mour}, $W$ is
$I$-cr and hence $W/W^{I} \hookrightarrow W$ obtained as an
$I$-splitting. Note that $E_I(W) = E_R(W) = E(W)$ is semistable of
degree 0 and since $W/W^{I}$ is an $I$-direct summand of $W$ the
associated bundle $E_I(W/W^{I})$ is a direct summand of the degree 0
semistable vector bundle $E_I(W)$. 

\noindent
This implies that $E_I(W/W^{I})$ is also semistable of degree 0.

Composing the section $s_1$ and the $I$-map $W \lr W/W^{I}$ we have:
\[
\overline s_1: E_I \lr W \lr W/W^{I}
\]
and $\overline w \in Im(\overline s_1)$. This gives a non-zero {\it
  unstable} section of $E_I(W/W^{I})$ which contradicts the stability
  of the bundle $E_I$ by Lemma \ref{stab}. 

This contradicts the assumption that the orbit $R \cdot w$ is {\it
  not} closed and completes the proof of the Proposition. Q.E.D.

\brem The theme in this section fits in with the general theme of
Kempf-Luna in the char.0 case. In char.0 the polystable bundle $E$
comes from an representation of $\pi_1(X) \lr G$. Let $R$ be the Levi
of an admissible parabolic and $E_R$ be as in \S9.  Then $E_R$ is
stable. So the representation $\pi_1(X) \lr G$ which factors via $R$
is irreducible . Let $M$ be the Zariski closure of the image. Then the
inclusion $M \hookrightarrow R$ is irreducible in the following
natural sense of \cite{mour} and \cite{tens}: namely, there exists no
parabolic subgroup $P \subset R$ such that $M \hookrightarrow P$.

In this case the proof of Proposition \ref{monodr} now follows easily
by results of Kempf. We need to check that the orbit $R \cdot w$ is
closed. Now $M$ is a reductive subgroup of $R$ which fixes $w$ since
$\pi_1(X)$ fixes $w$ (by classical local constancy). If $R \cdot w$ is
{\it not closed} then $R$ possesses a non-trivial one-parameter
subgroup and since $M$ fixes $w$ there exists a Kempf parabolic $P$
such that $M \hookrightarrow P \hookrightarrow R$ contradicting
irreducibility of $M \subset R$. (cf. \cite[Cor 4.4,4.5]{k}) \erem

\brem The Proposition \ref{monodr} appears in \cite{rr} but only in
char.0.  In \cite{rr} there is an error in the proof of the second
half of their theorem. Here we give a different proof of this and this
works in the situation when the action is separable which in
particular takes care of char.0 as well.  \erem

\section{Extension to the flat closure}

Recall that the section $s_K'(x)$ extends along the base section
$x_A$, to give $s_A'(x) = w_A$. The aim of this section is to prove the
following key theorem.
\noindent 
\bth\label{flatext} The section $s_K'$, extends to a section $s_A'$ of
 $E_A(G_A/H_A')$. In other words, the structure group of $E_A$ can be
 reduced to $H_A'$; in particular, if $H_k'$ denotes the closed fibre
 of $H_A'$, then the structure group of $E_k$ can be reduced to
 $H_k'$.  \eeth

\subsection{\it Saturated monodromy groups and Local constancy}
\bprop\label{lc} Let $E$ be a {\it polystable} principal $G$-bundle on
$X$.  Let $W$ be a $G$-module of low separable index, $w \in W$ and
$H' = Stab(w)$. Let $Y = G/H'$ the $G$-subscheme of $W$ defined by the
reduced subgroup $H' \subset G$. If $s$ is a section of $E(W)$ such
that for some $x \in X$, the evaluation at $x$, namely $s(x) = w$ is
in $E(Y)_x$, then the entire image of $s$ lies in $E(Y)$. In fact we
have a reduction of structure group to a reductive {\bf saturated}
subgroup $R_w$ of $H'$ and in particular, the reduced $R_w$-bundle is
stable.\eprop
\noindent
{\it Proof.}  Since the $G$-bundle $E$ is assumed polystable, by Def
\ref{poly}, there is an admissible reduction to a parabolic subgroup
$Q \subset G$ and a further reduction of structure group $E_R$, to a
Levi subgroup $R \subset Q$ with $E_R$ actually {\bf stable}.
 
Note that since $R$ is a Levi of a parabolic in $G$, the maximal torus
of $G$ and $R$ are the same.

Further, being a Levi of a parabolic $R$ is a saturated subgroup of
$G$. Since the height of the representation $G \lr SL(W)$ is low, it
follows that $W$ as an $R$-module is also of low height (cf.
\cite{mour} pp 22). 

Thus, we can conclude that $W$ as an $R$-module is also of {\it low
  separable index}.

Consider the $R$-bundle $E_R$ and the $R$-module $W$. We are given a
section $s :X \lr E(W) = E_R(W)$ such that at $x \in X$ $s(x) = w$ is
the given vector in $W$ with $Stab_G(w) = H'$.

By Proposition \ref{monodr}, since $E_R$ is stable, the orbit $R \cdot w =
F^{cl}$ is a closed orbit. Since the action of $R$ on $W$ is separable
the isotropy, $R_w = Stab_R(w)$ is reduced, and we have an isomorphism
$R.w \simeq R/R_w$. Note further that $R_w$ is {\it saturated} and
reductive. 

As one has observed in the previous proof the section $s$ takes its
values in the fibre $F$ and since $w \in F^{cl}$ we have the
following:
\[
s: E_R \lr F^{cl} \simeq R/R_w.
\]
\noindent
This gives a reduction of structure group of $E_R$ to $R_w$. We thus
have the following inclusion of bundles:
\[
E_{R_w} \hookrightarrow E_R \hookrightarrow E
\]
Note that $R_w =Stab_R(w) \subset Stab_G(w) = H'$. This inclusion
gives the required reduction of structure group of $E$ to $H'$ which
indeed comes as an extension of structure group from $E_{R_w}$.
Furthermore, $R_w$ is saturated and reductive. This complete the proof
of the Proposition.  Q.E.D.

\subsection{\it Completion of proof of Theorem \ref{flatext}}
\noindent
By Lemma \ref{chev} ,we have
\[
E_A(G_A/H_A') \hookrightarrow E_A(W_A).
\]  

The given section $s_K'$ of $E_K(G_K/H_K')$ therefore gives a section
 $u_K$ of $E(W_K)$.  Further, $u_K(x)$, the restriction of $u_K$ to $x
 \times T^*$, extends to give a section $u_A (x)$ of $E_x (W_A)$
 (restriction of $E_A(W_A)$ to $x \times T$).Thus, by Proposition
 \ref{dm}, and by the semistability
 of $E_l(W_A)$, the section $u_K$ extends to give a section $u_A$ of
 $E(W_A)$ over $X \times T$. 

Now, to prove the Theorem \ref{flatext} , we need to make sure that:
$$ 
\begin{array}{l}
\mbox{The image of this extended section $u_A$ actually lands in
$E_A(G_A/H_A')$} \\
\end{array} . \eqno (*)
$$
This would then define $s_A'$.

To prove ($\ast$), it suffices to show that $u_A(X \times l)$ lies in
$E_A(G_A/H_A')_l$ (the restriction of $E_A(G_A/H_A')$ to $X \times
l$).
 
Observe that, $u_A(x \times l)$ lies in $E_A(G_A/H_A')_l$ since
$u_A(x) = s_A'(x) = w_A$.
 
Observe further that, if $E_l$ denotes the principal $G$-bundle on
$X$, which is the restriction of the $G_A$-bundle $E_A$ on $X \times
T$ to $X \times l$, then $E_A(W_A)_l = E_A(W_A)|{X \times l}$, and we
also have

\[
\begin{array}{ccc}
E_A(G_A/H_A')_l &
\stackrel{\simeq}{\lr} &
E_l(G_k/H_k') \\
\Big\downarrow & &
\Big\downarrow \\
E_A(W_A)_l   & \stackrel{\simeq}{\lr} &
E_l(W)
\end{array}
\]
and the vertical maps are inclusions:
\[
E_A(G_A/H_A')_l \hookrightarrow E_A(W_A)_l, E_l(G_k/H_k')
\hookrightarrow E_l(W)
\]
where $E_l(W) = E_l \times^{H_k'} W$ with fibre as the $G$-module $W =
W_A \otimes k$.  Note that $G/H_k'$ is a $G$-subscheme $Y$ of $W$.  

\noindent
Recall that $E_l$ is polystable of degree zero.  Then, from the
foregoing discussion, the assertion that $u_A(X \times l)$ lies in
$E_A(G_A/H_A')$, is a consequence of Proposition \ref{lc} applied to
$E_l$. (Note that the group $H_k' = Stab_{G_k}(w_k)$ satisfies the
hypothesis of Proposition \ref{lc}).
 
Thus we get a section $s_A'$ of $E_A(G_A/H_A')$ on $X \times T$, which
extends the section $s_K'$ of $E_A(G_A/H_A')$ on $X \times T^*$.
This gives a reduction of structure group of the $G_A$-bundle $E_A$ on
$X \times T$ to the subgroup scheme $H_A'$ and this extends the given
bundle $E_K$ to the subgroup scheme $H_A'$.

In summary, we have extended the original $H_K$-bundle upto
isomorphism to a $H_A'$-bundle. The extended $H_A'$-bundle has the
further property that the limiting bundle $E_l'$ which is an
$H'_k$-bundle comes with a reduction of structure group to a reductive
and {\it saturated} subgroup $R_w$ of $H'_k$. Q.E.D

\brem The proof of Theorem \ref{flatext} is not as simple as in the
proof of Proposition \ref{dm}, since
\[ 
E_A(G_A/H_A') \hookrightarrow E_A(W_A)
\]
is {\it not} a closed immersion. The group scheme $H_A'$ is not
 reductive and therefore, we are {\it not} given a {\it closed}
 $G$-embedding of $G_A/H_A'$ in $G_A$-module $W_A$ (cf. Remark
 \ref{bull}).
\erem

\brem The reductive saturated subgroup $R_w$ plays the role of
``monodromy'' subgroup of the polystable $G$-bundle $E$. (cf.
\cite{bals}) 
\erem

\section{Potential good reduction}
Recall that by virtue of the separability of the action the group
scheme $H_A'$ is {\it smooth}. 

To complete the proof of the Theorem \ref{ssred}, we need to extend
the $H_K$-bundle to an $H_A$-bundle where $H_A$ denotes the reductive
group scheme $H \times \spec A$ over $A$.

\bprop\label{serre} There exists a finite extension $L/K$ with the
 following property: If $B$ is the integral closure of $A$ in $L$, and
 if $H_B'$ are the pull-back group schemes, then we have a morphism of
 $B$-group schemes
\[
\chi_B: H_B'~\longrightarrow H_B
\]
 which extends the isomorphism $\chi_L: H_L' \cong H_L$.
 \eprop
\noindent
 {\it Proof.} 
Observe first that the {\em lattice} $H_A'(A)$ is a {\it
 bounded subgroup} of $H_A(K)$, in the sense of the Bruhat-Tits
 theory \cite{bt}. Here, we make the identifications:
\[ 
H_K' \cong H_K~ \mbox {as K-group schemes}
\]
 Hence,
\[
H_A'(A) \subset H_K'(K) \cong H_K(K) = H_A(K)
\]
Then we use the following crucial fact: 
$$\left \{
\begin{array}{l}
\mbox {There exists a finite extension $L/K$ and an element $g \in
H_A'(L)$ such that}\\
\mbox{$g.H_A'(A).g^{-1} \hookrightarrow H_A(B)$.
}\\
\end{array} \right \} \eqno (*)
$$
This assertion is a consequence of the following result from,
(\cite{ser} Prop 8, p 546) (cf. also \cite {gille} Lemma I.1.3.2, or
\cite{lar} Lemma 2.4 ).

(Serre) There exists a totally ramified extension $L/K$ having the
following property: For every bounded subgroup $M$ of $H(K)$, there
exists $g \in H(K)$ such that $g.M.g^{-1}$ has {\it good reduction} in
$H(L)$ (i.e $h.M.h^{-1} \subset H(B)$, where $B$ is the integral
closure of $A$ in $L$).

For the sake of clarity we gather all the identifications of the
subgroups under consideration:
\[
H_A'(K) = H_K'(K)~\mbox {and}~H_A'(L) = H_B'(L) = H_L'(L)
\]
\[
H_A'(A) \subset H_B'(B)
\]
\[
H_A(B) = H_B(B) 
\]
Thus, we see that the isomorphism $\chi_L : H_L' \longrightarrow H_L$,
 given by {\it conjugation by} $g$, induces a map $\chi_L(B) : H_A'(A)
 \longrightarrow H_B(B)$. The crucial property to note is the
 following one:

 Given a rational point $\xi_k \in H_k'(k)$, there exists a point
 $\xi_A \in H_A'(A)$, and hence in $H_B'(B)$, which extends $\xi_k$,
 since $H_A'$ is smooth over $A$ and $k$ is algebraically closed.

 The proposition will follow by the following Lemma.  Let $A$, $B$ etc
 be as above.  

\blem Let $A$ be a complete discrete valuation ring with quotient
 field $K$.  Let $Z_A$ and $Y_A$ be $A$-schemes with $Z_A$ smooth. Let
 $\chi_L : Z_L \longrightarrow Y_L$ be a $L$-morphism such that
 $\chi_L(B) : Z_A(A) \longrightarrow Y_B(B)$. Then, the $L$-morphism
 $\chi_K$ extends to a $B$-morphism $\chi_B: U_B \longrightarrow Y_B$,
 where $U_B$ is an open dense subscheme of $Z_B$ which intersects all
 the irreducible components of the closed fibre $Z_k$.

In particular, if $Z_A$ and $Y_A$ are smooth and separated group
 schemes and if $\chi_L$ is a morphism of $L$-group schemes then
 there exists an extension $\chi_B : Z_B \longrightarrow Y_B$ as a
 morphism of $B$-group schemes.  
\elem
\noindent
 {\it Proof.} Consider the graph of $\chi_L$ and denote its schematic
 closure in $Z_B \times_B Y_B$ by $\Gamma_B$. Let $p : \Gamma_B
 \longrightarrow Z_B$ be the first projection. Then $p$ is an
 isomorphism on generic fibres. So, it is enough if we prove that $p$
 is invertible on an open dense $B$-subscheme $U_B$ of $Z_B$, which
 intersects all the components $C$, of the closed fibre $Z_k$. 

 We claim that, the map $p_k : \Gamma_k \longrightarrow Z_k$ is
 surjective onto the subset of $k$-rational points of each components,
 and this will imply that $p_k$ is surjective since $k$ is
 algebraically closed.  Note that $Z_A$ is assumed to be smooth and
 so, the closed fibre is reduced and also $k$ is algebraically
 closed. Thus, each $z_k \in Z_k(k)$ lifts to a point $z \in Z_A(A)
 \subset Z_B(B)$, $A$, being a complete discrete valuation ring. Since
 $\chi_L(B)$ maps $Z_A(A) \lr Y_B(B)$, we see that, there exists a $y
 \in Y_B(B)$ such that $(z,y) \in \Gamma_B(B)$. Thus, $z_k$ lies in
 the image of $p_k$. This proves the claim.

 In particular, by the well-known result of Chevalley on images of
 morphisms, the generic points, $\alpha$'s, of all the components $C$
 of $Z_k$ , lie in the image of $p_k$. Let $p_k(\xi) =
 \alpha$. Consider the local rings ${\cal O}_{\Gamma_{B},\xi}$ and
 ${\cal O}_{Z_{B},\alpha}$. Then by the above claim, the local ring
 ${\cal O}_{\Gamma_{B},\xi}$ dominates ${\cal
 O}_{Z_{B},\alpha}$. Since $Z_B$ is smooth and hence normal, for every
 $\alpha$ the local rings, ${\cal O}_{Z_{B},\alpha}$ are all discrete
 valuation ring's. Further, since $\Gamma_B$ is the schematic closure
 of $\Gamma_L$, it implies that $\Gamma_B$ is $B$-flat and $\Gamma_L$
 is open dense in $\Gamma_B$. Moreover, since $p$ is an isomorphism on
 generic fibres both local rings have the same quotient
 rings. Finally, since ${\cal O}_{Z_{B},\alpha}$ is a discrete
 valuation ring, we have an isomorphism of local rings. Therefore
 since the schemes are of finite type over $B$, we have open subsets
 $V_{i,B}$ and $U_{i,B}$ for each component of $Z_k$, which we index
 by $i$, such that $p$ induces an isomorphism between $V_{i,B}$ and
 $U_{i,B}$. This gives an extension of $\chi$ to open subsets
 $U_{i,B}$ for every $i$, with the property that these maps agree on
 the generic fibre. Since $Y_B$ is separated these extensions glue
 together to give an extension $\chi_B$ on an open subset, which we
 denote by $U_B$; this open subset will of course intersect all the
 components of the closed fibres of $Z_k$.
 
 The second part of the lemma follows immediately, if $Y_A$ is affine
 (which is our case). More generally, we appeal to the general theorem
 of A.Weil on morphisms into group schemes, which says that if a
 rational map $\psi_B$ is defined in codimension $\leq 1$ and if the
 target space is a group scheme then it extends to a global morphism.
 (cf. for example \cite{neron} pp 109). As we have checked above this
 holds in our case and this implies that as a morphism of schemes,
 $\psi_L$ extends to give $\psi_B: Z_B \lr Y_B$.

Further, by assumption $\chi_L$ is already a morphism of $L$-group
schemes and hence it is easy to see that the extension $\chi_B$ is
also a morphism of $B$-group schemes. This concludes the proof of the
lemma.

\brem  
Larsen in (\cite{lar}, (2.7) p 619), concludes from $(\ast)$, in the
$l$-adic case the statement of Proposition \ref{serre}. However, we
give a complete proof.
\erem

\brem\label{neron} In this section, by the assumption on the
separability index of the affine $(G,H)$-module we were able to conclude
that the flat closure $H_A'$ is indeed smooth. We observe that since
we are over char.$p$, in general the limiting fibre of the flat
closure $H_A'$ need not be {\it reduced}. This, as one knows is true
in char.0 by virtue of Cartier's theorem. Indeed more generally, given
a flat group scheme $H_A'$ with smooth generic fibre $H_K'$, there is
a construction due to Raynaud of what he calls the Neron-smoothening
of $H_A'$.  This exists as a smooth group scheme $H_A''$ with generic
fibre $H_K'' \simeq H_K'$ with the following universal property: given
any smooth $A$-scheme $D_A$ and an $A$-morphism $D_A \lr H_A'$, this
map factors uniquely via an $A$-morphism $D_A \lr H_A''$. In
particular, $H_A''(A) = H_A'(A)$.  Thus more generally without any
separability index assumptions, the proof of Proposition \ref{serre}
gives a morphism $H_B'' \lr H_B$. One is unable to make use of this
since the principal bundle $E_A'$ has structure group $H_A'$ and there
is no natural reason for its lifting to a principal $H_A''$-bundle.
(see \cite{neron}).  \erem

\section{Semistable reduction theorem}

Let {H} be a semi-simple algebraic group over $k$ an algebraically
closed field of char. p. Let $H \subset G = SL(V)$, be the
representation we have fixed in \S1. We retain all the notations of
\S7. The aim of this section is to prove the following theorem.
 
\bth\label{ssred} {(\em Semistable reduction)} Let $W$ be a finite
dimensional {\em affine $(G,H)$-module} associated to $H$ and $G$ and let
$p > \psi_G(W)$.  Let $H_K$ denote the group scheme $H \times \spec
K$, and $P_K$ be a semistable $H_K$-bundle on $X_K$.  Then there
exists a finite extension $L/K$, with $B$ as the integral closure of
$A$ in $L$ such that the bundle $P_K$, after base change to $\spec B$,
extends to a semistable $H_B$-bundle $P_B$ on $X_B$.  \eeth

\noindent
{\it Proof.}  First by Proposition \ref{flatext} we have an
$H_A'$-bundle which extends the $H_K$-bundle upto isomorphism. Then,
by Proposition \ref{serre}, by going to the extension $L/K$ we have a
morphism of $B$-group schemes $\chi_B : H_B' \longrightarrow H_B$
which is an isomorphism over $L$. Therefore, one can extend the
structure group of the bundle $E_B'$ to obtain an $H_B$-bundle $E_B$
which extends the $H_K$-bundle $E_K$.

We need only prove that the fibre of $E_B$ over the closed point is
indeed {\it semistable}. This is precisely the content of Proposition
\ref{lastpush} below. Q.E.D.

\brem We remark that this is fairly straightforward in char.0 since it
comes as the extension of structure group of $E_l'$ by the map $\chi_k
: H_k' \lr H_k$. We note that in char.0, $E_l'$ is the $H_k'$-bundle
obtained as the reduction of structure group of the polystable vector
bundle $E(V_A)_l$ and so remains semistable by any associated
construction (cf. Proposition 2.6 of \cite{bals}). In our situation
this becomes much more complex and we isolate it in the following
proposition. \erem

\bprop\label{lastpush} The limiting bundle, namely the fibre of $E_B$
over the closed point is semistable.
\eprop
\noindent
{\it Proof.} Recall from Proposition \ref{lc}, that the limiting
bundle of the family $E_B'$ namely $E_l'$, had the property that it
had a further reduction of structure group to a reductive and
saturated group $R_w$ of $H_k'$ and hence of $G_k = G$. Thus the
representation $R_w \lr G_k$ is also {\it low height} by (\cite{mour}
pp 25 ).  Further, by the low height property, the representation $R_w
\lr G = SL(V)$ is completely reducible. 
 
Observe further that since $H_k'$ is {\em not reductive} (cf. Remark
\ref{bull} above), there exists a proper parabolic subgroup $P \subset
G_k$ such that $H_k' \subset P$. This follows by the theorem of
Morozov-Borel-Tits (cf. \cite{bort}). Therefore the subgroup $R_w
\subset H_k' \subset P$. Now since $R_w \lr G$ is completely
reducible, $R_w \subset P$ implies that $R_w \subset L$ for a Levi
subgroup $L \subset P$.

Now $R_w$ is a saturated reductive subgroup of $G$. Therefore, since
$p > \psi_G(W)$, by Lemma \ref{heightlie} (and Remark
\ref{heightlie1}) below we see that the modules $Lie G_k$ and $Lie
H_k'$ are low height modules for $R_w$ and in particular completely
reducible.

Now $R_w$ is a saturated group and the connected component of
identity, $R_w^{0}$, is reductive by Proposition \ref{lc}.  Since
$R_w$ is saturated as a subgroup of $G$ by height considerations, the
modules $Lie G_k$ and $Lie H_k'$ are low height modules for $R_w^{0}$
as well. (cf. Remark \ref{saturate})

Thus by Remark \ref{cline}, we have the following:
\[
H^{i}(R_w^{0},Lie(H_k')) = H^{i}(R_w^{0},Lie(G_k)) = 0
\]
for all $i \geq 1$.

Recall that by Remark \ref{saturate} the {\it saturatedness} of $R_w$
implies that the index $[R_w : R_w^{0}]$ is prime to the
characteristic $p$. 

Therefore if we denote $R_w/R_w^{0}$ by $I_w$, we see that the order
of $I_w$ is prime to $p$. Hence we have the following vanishing of
cohomology:
\[
H^{i}(I_w,Lie(H_k')) = H^{i}(I_w,Lie(G_k)) = 0
\]
for all $i \geq 1$. (For this classical result cf. \cite{cartan}
p. 237.)

Putting together the above results, we can conclude the following:
\[
H^{i}(R_w,Lie(H_k')) = H^{i}(R_w,Lie(G_k)) = 0
\]
for all $i \geq 1$.

This implies, by the infinitesimal lifting property of (\cite{SGA} Exp.III
Cor 2.8) that if we consider the product group scheme $R_{w,B} = R_w \times
Spec(B)$, then the inclusion 
\[
i_k: R_w \hookrightarrow H_k' \hookrightarrow G_k
\]
lifts to an inclusion

\[
i_B: R_{w,B} \hookrightarrow H_B' \hookrightarrow G_B
\]
where the generic inclusion is defined upto conjugation by the
inclusion over the residue field. 

Denote the above composite by:
\[
i_{1,B}: R_{w,B} \hookrightarrow G_B
\]
By Proposition \ref{serre}, we also have a morphism $\chi_B: H_B' \lr
H_B$, which is an isomorphism over the function field $L$. We have the
following diagram:
\begin{center}
\begin{picture}(400,325)(0,20)
\put(25,290){\ctext{$R_{w,B}$}}
\put(75,290){\vector(1,0){200}}
\put(325,290){\ctext{$H_B'$}}
\put(325,250){\vector(0,-1){200}}
\put(75,265){\vector(1,-1){200}}
\put(325,25){\ctext{$H_B$}}
\put(125,150){\ctext{$j_B$}}
\put(350,150){\ctext{$\chi_B$}}
%\put(250,200){\ctext{$\circledarrowright$}}
\put(175,320){\ctext{$i_B$}}
\end{picture}
\end{center}

We note that we also have an inclusion $H_B \hookrightarrow G_B$
coming from the original representation $H \hookrightarrow G$. In other
words we have another morphism
\[
j_{1,B}: R_{w,B} \lr G_B
\]
Thus, we get the following diagram:

\begin{center}
\begin{picture}(400,325)(0,20)
\put(25,290){\ctext{$R_{w,B}$}}
\put(75,290){\vector(1,0){200}}
\put(325,290){\ctext{$G_B$}}
%\put(325,250){\vector(0,-1){200}}
\put(75,265){\vector(1,-1){200}}
\put(325,25){\ctext{$G_B$}}
\put(125,150){\ctext{$j_{1,B}$}}
%\put(350,150){\ctext{$\chi_B$}}
%\put(250,200){\ctext{$\circledarrowright$}}
\put(175,320){\ctext{$i_{1,B}$}}
\end{picture}
\end{center}
(We remark that there is no vertical arrow to complete the above
diagram!) 

Note that over the function field $L$ the maps $j_{1,L}$ and $i_{1,L}$
coincide upto conjugation. Thus by the cohomology vanishing stated
above and the rigidity of maps (\cite {SGA} Exp III Cor 2.8), the maps
over the residue fields are also conjugates.

Consider the bundle $E_l'$ which comes equipped with a reduction to
$R_w$ and is semistable as an $R_w$-bundle (cf. Prop \ref{lc}). 

Since the representations $i_{1,k}:R_w \hookrightarrow G_k$ and
$j_{1,k}:R_w \hookrightarrow G_k$ are conjugate it follows that the
associated $G_k$-bundles $E_l'(j_{1,k})$ and $E_l'(i_{1,k})$ are
isomorphic. Therefore since $E_l'(i_{1,k})$ is semistable so is
$E_l'(j_{1,k})$. In particular, since the morphism $j_{1,k}:R_w
\hookrightarrow G_k$ factors via $H_k$, the associated $H_k$-bundle
$E_l'(j_k)$ is semistable.  This implies that the induced bundle $E_B$
is a family of semistable $H_B$-bundles. This completes the proof of
the Theorem \ref{ssred}. Q.E.D.

\brem\label{del} Let $H \subset G$, where $G$ is a linear group. In
the notation of \S2 let $F_H$ and $F_G$ stand for the functors
associated to families of semistable bundles of degree zero. (cf.
Proposition \ref{dm}). The inclusion of $H$ in $G$ induces a morphism
of functors $F_H \lr F_G$. We remark that the semistable reduction
theorem for principal $H$-bundles {\bf need not} imply that the
induced morphism $F_H \lr F_G$ is a proper morphism of
functors. Indeed, this does not seem to be the case. However, it does
imply that the associated morphism at the level of moduli spaces is
indeed proper (cf. Theorem \ref{moduli}).  \erem

\subsection{\it Some remarks on low height modules}

\blem\label{heightcomps} Let $H \subset G = SL(V)$ be a low height
representation. Let $W$ be a low height $G$-module such that $G/H$ is
embedded as a closed orbit in $W$ (cf. Def \ref{chevmod}). Suppose
that the subspace $V^{H} \subset V$ of $H$-fixed vector in $V$ is the
zero subspace. Then $W$ contains direct summand different from $V$ and
$V^*$. (Note that by the low height assumptions all modules are
completely reducible.)\elem
\noindent
{\it Proof.} For if $W = \oplus V$, then the vector $w \in W$ which
has a closed $G$-orbit and whose isotropy is $H$ projects onto a
vector $v \in V$ fixed by $H$. But by assumption, the subspace $V^{H}
= 0$. Hence $W$ cannot be a direct sum of copies of $V$. We also
observe that this implies $(V^{*})^{H} = 0$ as well and therefore $W$
is not the direct sum of $V^*$'s alone.

\blem\label{heightcomps1} Let $R \subset G = SL(V)$ be a reductive
saturated subgroup of $G$ that is contained in the Levi of a parabolic
subgroup of $G$. Let $W$ be a low height $G$-module that contains a
component not isomorphic to $V$ and $V^*$. Then $Lie(G)$ and $Lie(H)$
are low height $R$-modules, in particular completely reducible.  \elem
\noindent
{\it Proof.} Let $n = dim(V)$. Since $W$ contains a component other
than $V$ and $V^*$, $ht_G(W) \geq 2(n - 2)$.

Hence $W$ being a low height $G$-module we have $ p > 2(n - 2)$. Since
$R$ is not irreducible in that $R$ is contained in a certain Levi
subgroup $L \subset P$ of a parabolic subgroup $P \subset G$, it
follows that $ht_R(V) < ht_G(V) = n - 1$.

Hence $ht_R(V \otimes V^{*}) \leq 2(n - 2) < p$. In other words, $V
\otimes V^* = Lie(G)$ is a low height $R$-module. Note also that
$ht_R(Lie (H)) \leq ht_R(Lie(G))$. Q.E.D.

\blem\label{heightlie} Let $H \subset G = SL(V)$ be a low height
representation. Let $W$ be a low height $G$-module such that $G/H$ is
embedded as a closed orbit in $W$ (cf. Def \ref{chevmod}).  Let $R
\subset G = SL(V)$ be a reductive saturated subgroup of $G$ that is
contained in the Levi of a parabolic subgroup of $G$. Assume that
$V^{H} \subset V^{R}$. Then $Lie(H)$ and $Lie(G)$ are low height
$R$-modules.  \elem
\noindent
{\it Proof.} Let $V'$ be the subspace complementary to $V^{H}$ in $V$.
Let $n = dim(V)$ and $n' = dim(V')$. Let $G' = SL(V') \subset G$.
Then the representation $H \hra SL(V) = G$ factors through $H \hra
SL(V') = G'$. Moreover $G'$ is saturated in $G$ (being the semisimple
part of a Levi subgroup of a parabolic subgroup) and therefore $V$ and
$W$ are low height $G'$-modules (by Remark \ref{saturate}).

By the choice of $V'$, we have $(V')^{H}= 0$. Since $V^{H} \subset
V^{R}$ we see that $R \subset G'$. Therefore the $G'$-orbit gives a
closed embedding of $G'/H$ in $W$. It follows by Lemma
\ref{heightcomps} that $W$ contains summands other than $V'$ and
${V'}^*$.

Hence by Lemma \ref{heightcomps1} $Lie(G')$ and $Lie(H)$ are low
height $R$-modules.  Now the result follows because $ht_H(V) =
ht_H(V')$ and hence $ht_R(V) = ht_R(V')$.  This works for the duals as
well, i.e $ht_R(V^*) = ht_R({V'}^*)$. By additivity of heights we see
that
\[
ht_R(V \otimes V^*) = ht_R(V' \otimes {V'}^*) < p
\]
since $Lie(G')$ is a low height $R$-module. Q.E.D.

(cf. \cite{mour} p. 27 for some of the computations made here)

\brem\label{heightlie1} We note that the subgroup $R_w$ to which we
apply Lemma \ref{heightlie} satisfies the condition of the Lemma,
especially the condition that $V^{H} \subset V^{R_w}$. This follows
since $H_A'$ is the {\em flat closure} of $H_K'$ in $G_A$ and since
$R_w \subset H_k'$. In fact, for the purposes of Prop \ref{lastpush}
or the semistable reduction theorem one could have worked with $G' =
SL(V')$ instead of $G$. In that case it is clear that the flat closure
of $H_K'$ is actually realised in $G'_A$ itself.\erem

\subsection{\it Irreducibility of the moduli space}

 We first remark that the semistable reduction theorem Theorem
\ref{ssred} holds in fact in a slightly more general setting as well.

\bcor \label{gssred} Let $\cal X \lr S$ be a smooth family of curves
parametrised by $S = Spec A$ where $A$ is a complete discrete
valuation ring with char.$K$ = 0 and residue characteristic $p$.
Suppose further that $p > \psi(W)$ as in Theorem \ref{ssred}. Let
$H_S$ be a reductive group scheme obtained from a split Chevalley
group scheme $H_{\bf Z}$.  Suppose further that we are given a family
of semistable principal $H_K$-bundles $E_K$ over ${\cal X}_K$. Then,
there exists a finite cover $S' \lr S$ such that the family after
pull-back to $S'$ extends to a semistable family $E_{S'}$.  \ecor
\noindent
{\it Proof.} The proof of Theorem \ref{ssred} goes through with some
minor modifications.

We then have \bcor\label{irred} Let $H$ be simply connected. Then for
$p > \psi_G(W)$ the moduli spaces $M(H)$ of principal $H$-bundles is
irreducible.  \ecor
\noindent
{\it Proof.}  The proof of this is now standard once Cor \ref{gssred}
is given and one knows the fact over fields of char 0. The argument
very briefly runs as follows: The first point is to observe that given
the prime bounds, namely $p > \psi_G(W)$, the moduli scheme can be
constructed as in \S4 over $S = {\bf Z} - \{p \leq \psi_G(W)\}$. Call
this scheme $M(H)_S$. Then Cor \ref{gssred} implies that $M(H)_S$ is
projective and further, GIT (cf. \cite{ses1}) implies that the
canonical map $M(H) \lr M(H)_S \otimes k$ is a bijection on $k$-valued
points. Further, since $M(H)_S$ is projective and connected over the
generic fibre (by char 0 theory), Zariski's connectedness theorem
implies that the closed fibre $M(H)_S \otimes k$ is also connected and
hence so is $M(H)$. Now observe that the quot scheme $Q''$ constructed
in \S4 is easily seen to be {\it smooth} by some standard deformation
theory. Hence $M(H)$ is normal and connected and therefore {\it
irreducible}.

\small
\vspace{.01in}
\begin{flushleft}
V.Balaji \\
Chennai Mathematical Institute\\
92, G.N.Chetty Road, Chennai-600017 \\
INDIA \\
E-mail:
balaji@cmi.ac.in \\ [2mm]

A.J.Parameswaran \\
School of Mathematics,\\ 
Tata Institute of Fundamental Research.\\
Homi Bhabha Road,\\
Mumbai-40005 \\
INDIA \\ 
E-mail:
param@math.tifr.res.in \\

\end{flushleft}

\end{document}